\documentclass[10pt,oneside,onecolumn,final,leqno]{amsart}
\usepackage{amscd}
\usepackage{amssymb}
\usepackage{amsfonts}
\usepackage{pb-diagram}
\newtheorem{cor}{Corollary}[section]
\newtheorem{defin}{Definition}[section]
\newtheorem{teor}{Theorem}[section]
\newtheorem{prop}{Proposition}[section]
\newtheorem{lem}{Lemma}[section]
\newtheorem{oss}{Remark}[section]
\newcommand{\uno}{\boldsymbol{u}}
\newcommand{\uu}{(u_{2\alpha-1} , u_{2\alpha} )}
\newcommand{\uua}{(u_{2\alpha-1} , u_{2\alpha} )}
\newcommand{\un}{u_{2\alpha-1}}
\newcommand{\um}{ u_{2\alpha} }
\newcommand{\barun}{{\overline{u}}_{2\alpha-1}}
\newcommand{\barum}{{\overline{u}}_{2\alpha}}

\newcommand{\sq}{{}^{1\pm 2}\square_{1\pm 4}^{1\pm 3}}

\newcommand{\sqpm}{{}^{1\pm 2}\square_{1\pm 4}^{1\pm 3}}

\newcommand{\cost}{cos{\, \theta}}

\newcommand{\sit}{sin{\, \theta}}

\newcommand{\bfl}{\begin{flushleft}}
\newcommand{\efl}{\end{flushleft}}
\newcommand{\ka}{\mathrm{K\ddot{\mathrm{a}}hler}}

\newcommand{\grm}{Gr_4 (\mathbb{R}^8)}

\newcommand{\mal}{M_{\alpha} }
\newcommand{\zw}{(\underline{z}, \underline{w})}

\begin{document}

\title[Self Dual Einstein Orbifolds] {Self Dual Einstein Orbifolds with Few Symmetries\\
 as Quaternion K\"ahler Quotients}
\author{Luca Bisconti and  Paolo Piccinni}
\date{\today}

\dedicatory{To the memory of Krzysztof Galicki}

\address{L. B.:  I.N.G.V.\\ Sezione di Pisa\\
via della Faggiola, 1\\ 56126 Pisa\\ Italia.}
\email{bisconti@pi.ingv.it}
\address{P. P.:  Dipartimento di Matematica\\ Universit\`a di Roma 'La Sapienza'\\Piazzale Aldo Moro 2\\
I-00185 Roma\\ Italia.}
 \email{piccinni@mat.uniroma1.it}

\begin{abstract}
We construct a new family of compact orbifolds
$\mathcal{O}^4(\Theta)$ with a positive self dual
Einstein  metric 
and a one-dimensional
group of isometries.  Together with another family, introduced in \cite{bg} and here denoted by  $\mathcal{O}^4(\Omega)$, these examples classify all $4$-dimensional orbifolds that are quaternion  K$\ddot{\mathrm{a}}$hler quotients  by a torus of real Grassmannians.
\end{abstract}
\maketitle

\section{Introduction}
A  classical  theorem of Hitchin states that $S^4$ and $\mathbb CP^2$, with their symmetric metrics, are the only compact positive self dual Einstein (SDE) $4$-manifolds \cite{hi}. A classification of compact positive self dual Einstein $4$-orbifolds appears much harder, and at the present not fully understood. 
 
 First examples of compact self dual Einstein $4$-orbifolds of positive scalar curvature were constructed by Galicki and Lawson via their also now classical quaternion K\"ahler quotient construction. Presently, known examples of such orbifolds include: (i) the $SO(3)$-invariant, cohomogeneity one orbifold metrics on $S^4$ discovered by Hitchin \cite{hd}; (ii) the toric orbifold metrics constructed by Boyer, Galicki,  Mann and Rees as quaternion K\"ahler quotients of quaternionic projective spaces \cite{mr}; (iii) the $S^1$-invariant orbifold metrics of Galicki and Nitta \cite{gn}.
 
The toric orbifold metrics mentioned in (ii) include as special cases the Galicki-Lawson metrics on weighted complex projective spaces $\mathbb C P^2(p,q,q)$. All these toric metrics have been completely classified through quaternion $\ka$ quotients first by Bielawski \cite{bi} in a special case and more generally by Calderbank and Singer \cite{cs}. Compact positive SDE orbifolds with a one-dimensional isometry group are known just in a few cases, and the only known examples seem to be the ones mentioned in (iii) and the family constructed in \cite{bg} again through a quaternion $\ka$ quotient.

The present paper is devoted to constructing a new family of positive SDE metrics with a one-dimensional isometry group on compact orbifolds. We show that these new examples, together with the ones constructed in \cite{bg}, classify all such orbifolds that can be obtained as toric quotients from the quaternion $\ka$ Grassmannians $Gr_4 (\mathbb{R}^{n+1})\cong\frac{SO(n+1)}{SO(n-3)\times SO(4)}$.
In fact, actions by a $k$-dimensional torus sitting inside the maximal torus of $SO(n+1)$ and leading to a 4-dimensional positive SDE quotient orbifold give necessarily $n+1 = 6,7,8$ (cf. Section 2). Indeed the first case $n+1 = 6$ gives rise to SDE orbifolds with $T^2-$ symmetries \cite{mr}. The remaining cases $n+1 = 7, 8$ give rise to intermediate circle quotients related to the groups $G_2$ and $Spin(7)$, respectively: \cite{sk}, \cite{lp}.
 
Let  $\mathcal{M}^{4n}$ be a quaternion $\ka$ manifold of positive scalar
 curvature. Recall the diagram of fibrations, all consistent with the respective
quotient constructions:
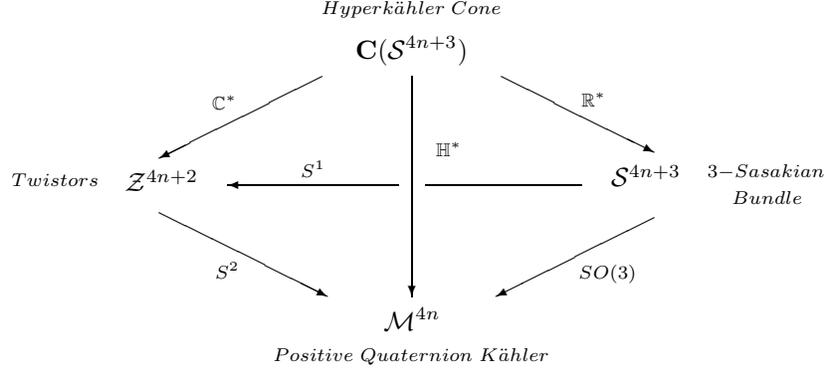
\[
   \dgARROWLENGTH=0.05\dgARROWLENGTH
   \begin{diagram}
   \node[2]{\overset{ \left.\begin{array}{c} {}^{ Hyperk\ddot{a}hler \; Cone } \\ \end{array}\right.}{\mathbf{C}(\mathcal{S}^{4n+3})}}\arrow[2]{s,r}{ \left.\begin{array}{c}\overset{\mathbb{H}^*}{} \\ \\ \\ \end{array}\right.}
      \arrow{sw,t}{ \mathbb{C}^*}  \arrow{se,t}{\mathbb{R}^*} \\
   \node{\left. \begin{array}{cc} \overset{Twistors}{} & \mathcal{Z}^{4n+2}\\ \end{array} \right.  } \arrow{se,b}{S^2}  \node{}
   \arrow{w,t}{S^1}     \node{\left. \begin{array}{cc}\mathcal{S}^{4n+3} &  \underset{{Bundle} }{{}^{3-Sasakian}} \end{array} \right.} \arrow{w,-} \arrow{sw,b}{SO(3)} \\
                           \node[2]{\underset{ \left.\begin{array}{c}  {}_{Positive \; Quaternion \; K\ddot{a}hler} \\ \end{array} \right. }{\mathcal{M}^{4n} } }   
   \end{diagram}
\]
In particular, one can consider weighted action of tori
$T^{k}_{\Theta}\subset SO(n+1)\subset$
$Sp(n+1)$ on the 3-Sasakian sphere $ S^{4n+3}\subset\mathbb{H}^{n+1}$, $k=[\frac{n+1}{2}]-1$,
where $\Theta$ is a $(k-1) \times k$ integral matrix. In Section 2
we describe the details when $n=7$ and the acting torus is $3-$dimensional. This gives the following:

\vspace{0,3 cm}

\noindent  $\mathbf{Theorem\,\, A}$   \emph{
 Let $\Theta\in M_{3\times 4}(\mathbb{Z})$
  be an integral matrix such that each of its  $3\times 3$ minor determinants $\Delta_{\alpha\beta\gamma}$
does not vanish. Moreover, assume that their sum is non zero, that none of them is equal to the sum of the other three, and that none of the sums of two of them is  equal to the sum of the other two. Then, for each such a matrix $\Theta$,
there exists a  compact self dual Einstein  $4-$dimensional
 orbifold $\mathcal{O}^4(\Theta )$ with positive scalar curvature and a
 one-dimensional group of isometries.$\square$}

\vspace{0,3 cm}

In fact, the analysis of the action shows that
 no choice of matrix $\Theta$ gives a smooth quaternion $\ka$ quotient of the corresponding quaternion $\ka$ base $\mathbb HP^7$, and this is the case also for $3-$Sasakian quotient metric on the mentioned $SO(3)$-bundle. The  singularities of  the quotient $\mathcal{O}^4(\Theta)$
can in fact be more conveniently described through the singular locus on its twistor space $\mathcal{Z}^6(\Theta)$. To give a proper formulation of this, denote by $\widetilde{G}$ the group $Sp(1) \times T^{3}_{\Theta} \times U(1)$. It acts in a natural way on the quaternionic vector space $\mathbb{H}^{8} $, and denote by $u_\alpha = z_\alpha +j w_\alpha (\alpha= 1,...,8)$ its coordinates. In Section 3 we prove the following:

\vspace{0,3 cm}

\noindent $\mathbf{Theorem\,\, B}$  \emph{Let $\mathcal{Z}^6(\Theta)$ be the twistor space of the self dual Einstein orbifold $\mathcal{O}^4(\Theta)$. Then the singular locus $\Sigma(\Theta)$ of $\mathcal{Z}^6(\Theta)$ contains at most the following sets.}
\begin{itemize}
\item[(i)] \emph{\noindent Two spheres $S^2$, whose isotropy
depends only on an algebraic sum of the determinants $\Delta_{\alpha\beta\gamma}$.} \\
\item[(ii)]  \emph{\noindent Further $22$ disjoint spheres $S^2$, whose isotropy depends on the minor determinants $\Delta_{\alpha\beta\gamma}$. These $2$-spheres are obtained as $\widetilde{G}-$ quotients from strata  ${S}{}^{\alpha \beta \gamma}_{\delta}$, ${S}{}^{\alpha}_{ \beta \gamma \delta}$ or ${S}{}^{\alpha \beta }_{\gamma \delta}$ on $\mathbb{H}^8$, that are $\widetilde{G}$ orbits of loci where some pairs of complex coordinates $(z_{2\alpha -1}, z_{2\alpha})$ or $(w_{2\alpha -1}, w_{2\alpha})$ are zero. For example:}
\begin{align*}
\qquad \qquad &{S}{}^{123}_{4} =\widetilde{G}\cdot \Bigg\{
\left( \begin{array}{cc|cc|cc|cc}
z_{1} & z_{2}  & z_{3} & z_{4}   &    z_{5} & z_{6}        & 0    &  0   \\
0   &    0 &  0  & 0     &   0   &  0          & w_{7}  &  w_{8}\\
\end{array}\right) \Bigg\},  \\
\qquad \qquad &{S}{}^{1}_{234} =\widetilde{G}\cdot \Bigg\{
\left( \begin{array}{cc|cc|cc|cc}
z_{1} & z_{2}  &   0  & 0     &   0   &  0  & 0    &  0   \\
0   &    0 &   w_{3} & w_{4}   &    w_{5} &w_{6}             & w_{7}  &  w_{8}\\
\end{array}\right) \Bigg\},  \\
\qquad \qquad &   {S}{}^{12}_{34} =\widetilde{G}\cdot \Bigg\{
\left( \begin{array}{cc|cc|cc|cc}
z_{1} & z_{2}  & z_{3} & z_{4}       & 0        &  0           & 0    &0\\
0     &    0   &  0    & 0           &    w_{5} & w_{6}     & w_{7}  &  w_{8}\\
\end{array}\right) \Bigg\}.
\end{align*} 
 \emph{\noindent One obtains in this way eight strata ${S}{}^{\alpha \beta \gamma}_{\delta}$ and ${S}{}^{\alpha}_{ \beta \gamma \delta}$ each of which intersects the zero set of the moment map in two connected components, and six further strata ${S}{}^{\alpha \beta }_{\gamma \delta}$ having a connected intersection with the zero set. All quotients by $\widetilde{G}$ of  these connected components are spheres $S^2$.} \\

\item[(iii)] \emph{\noindent  Three sets of  at most four points. The points of each set are joined by one of the  $2-$spheres
${S}{}^{\alpha\beta}_{\alpha\gamma}\cap N(\Theta)\big/\widetilde{G}$, where $N(\Theta) \subset S^{31} \subset \mathbb H^8$ is the zero set of the moment map associated to the $Sp(1) \times T^3_\Theta$- action.}\\
\end{itemize}

\vspace{0.3 cm}

We describe also singularities for 
quotients appeared in \cite{bg}. We use here the notation $\Omega \in M_{2 \times 3}(\mathbb Z)$ for the matrix of weights, entering in the $\tilde{G}^{\Omega}= T^2_{\Omega}\times
Sp(1)\times U(1)$ action on the quaternionic vector space $\mathbb H^7$ (cf. Section 4). We denote by $\Delta_{\alpha\beta}$, $\alpha, \beta = 1,2,3$,
the minor
determinants of the matrix $\Omega$. It is proved in \cite{bg} that, under some hypotheses on the $\Delta_{\alpha\beta}$, a positive SDE orbifold $\mathcal{O}^4(\Omega)$ with a one-dimensional group of isometries can be constructed. In Section 4 we describe the singular locus at twistor level by proving the following:

\vspace{0.3 cm}
$\mathbf{Theorem\,\, C}$ \emph{ Let
 $\mathcal{Z}^6(\Omega)$ be the twistor space of the self dual Einstein orbifold
  $\mathcal{O}^4(\Omega)$. The singular locus 
$\Sigma(\Omega)\subset \mathcal{Z}^6(\Omega)$
  contains at most the following sets: } 
\begin{itemize}
\item[(i)]  \emph{one sphere $S^2$, whose isotropy
depends only on one of the possible algebraic sums of the  $\Delta_{\alpha\beta}$};   
\item[(ii)]  \emph{$12$ points, coming from
  the following strata of the action of 
$\tilde{G}^{\Omega} $ on $N(\Omega)$ on $\mathbb{H}^7$}:   
 \[ 
 S^{12}_3 =\widetilde{G}^{\Omega}\cdot \Bigg\{
\left( \begin{array}{c|cc|cc|cc}
0   & z_2  & z_3   & z_4   &    z_5 &  0    &  0   \\
0   &    0 &  0  & 0    &   0    & w_6    &  w_7\\
\end{array}\right) \Bigg\},  
\]
\emph{and the similarly defined} $S^{13}_2,S^{23}_1, S^{1}_{23},S^{2}_{13}, S^{3}_{12}.   $

\emph{Any stratum ${S}{}^{\alpha  \beta }_{\gamma}$ or ${S}{}^{\alpha}_{ \beta \gamma}$ intersects
 the zero set of the moment map in
 two connected components. Each of these connected components
gives rise to a singular point at the twistor level.
   Moreover, for each of these points the isotropy depends only on
one of the minor determinants $\pm \Delta_{\alpha \beta} $.}
\end{itemize}

\vspace{0.3 cm}

When some of the minor determinants
  $\Delta_{\alpha\beta\gamma}$ or of their algebraic sums $\sq$ are $\pm 1$ in Theorem B, or similarly when some of the $\Delta_{\alpha\beta}$ or their algebraic sums are $\pm 1$ in Theorem C, then  the singular loci
 $\Sigma(\Theta)$  and  $\Sigma(\Omega)$ do not contain the corresponding sets.

 The
 comparison between singularities in the two cases shows that Theorem A 
gives rise to a new family of positive SDE orbifolds with a
 one-dimensional
 group of isometries.
\vspace{0.3 cm}

\noindent 
\emph{Acknowledgement and Provenance.} This paper is based on the first author's doctoral thesis, defended at Roma Tor Vergata in 2007 \cite{bis}. Both authors express their gratitude to Krzysztof Galicki, for his decisive encouragement and many helpful suggestions and discussions. This paper is dedicated to his memory.

\vspace{1 cm}

\section{The Quotient Orbifolds $\mathcal{O}^4 (\Omega)$ and $\mathcal{O}^4 (\Theta)$}

A family of $4-$dimensional positive 
SDE orbifolds with
one-dimensional group of isometries has been constructed in \cite{bg}. We denote here by  $\mathcal{O}^4(\Omega)$ these orbifolds, a notation that allows to distinguish them from the new orbifolds $\mathcal{O}^4(\Theta)$ that will be introduced in the present paper. We recall that the  
$\mathcal{O}^4(\Omega)$ are quaternion $\ka$ quotients, via a
$Sp(1)\times T^2_\Omega$ action with convenient weight matrix $\Omega \in M_{2 \times 3}(\mathbb Z)$ on the torus factor,
of the quaternionic projective space $\mathbb HP^6$. An alternative quotient construction of the same $\mathcal{O}^4(\Omega)$
is through the action of the weighted
$2-$torus $T^2_\Omega$  on the oriented
Grassmannian $Gr_4 (\mathbb{R}^{7})$. Under suitable assumptions for the weight matrix
$\Omega$, orbifold quotients $\mathcal{O}^4(\Omega)$ are obtained.

A similar construction by a torus action can be carried out on any quaternion-K$\ddot{\mathrm{a}}$hler oriented
Grassmannian

\begin{equation}         \label{C0 : 12} 
Gr_4 (\mathbb{R}^{n+1})\cong\frac{SO(n+1)}{SO(n-3)\times SO(4)},
\end{equation}
and if a $4$-dimensional quotient is desired, one has to look at actions of  $(n-4)-$
dimensional
tori. The dimension $[\frac{n+1}{2}]$ of the maximal torus in $SO(n+1)$ shows that the possibility of introducing weights in the torus action yields the inequality $ n-4 < [\frac{n+1}{2}].$
Thus:
 
\begin{equation}                       \label{C0 : 14}
 \left. \begin{array}{l}
           {}\,\,\, n < 9  \,\,\,\,(n+1\,\,\mathrm{even})
            \end{array} \right.
  \quad \quad \,\, \mathrm{or} \quad \quad \quad \quad \left. \begin{array}{l}
           {}\,\, n < 8  \,\,\,\,  (n+1\,\, \mathrm{odd}),
         \end{array} \right.
\end{equation}
so that the only Grassmannians that can admit such quotients are:

\begin{equation}  \label{C0: 14}
 \left. \begin{array}{lllll}
   Gr_4 (\mathbb{R}^6)\cong Gr_2(\mathbb{C}^4),  &  {}  \quad  {} & Gr_4 (\mathbb{R}^7),
       &    {}  \hspace{1 cm} \mathrm{and} \quad {}& Gr_4 (\mathbb{R}^8). \\
       \end{array} \right.
\end{equation}
The first two cases in (\ref{C0: 14}) have been examined
in \cite{mr}  and in  \cite{bg}, respectively.  
The present paper is 
devoted to the third case and to its comparison with the second one (circle quotients of $Gr_4 (\mathbb{R}^6)$ have a two-dimensional group of isometries, and are therefore a priory distinct from orbifolds in the other two families). Thus our main choice is
the Grassmannian $\grm$, acted on
by a $3-$torus $T^3_{\Theta}\subset T^4 \subset SO(8)$,
where $\Theta$ is a $3\times 4$ integral weight matrix. 

The action of $T^3_{\Theta}$ is conveniently described
through  $2\times 2$
  block diagonal
 matrices:   

 \begin{equation}         \label{C0 : 5}
 A(\Theta)=\left(\begin{array}{c|c|c|c}
 A(\theta_1) & 0 & 0 & 0  \\
 \hline
  0 & A(\theta_2) & 0 & 0  \\
  \hline
 0  & 0   & A(\theta_3)  & 0\\
  \hline
 0 & 0 & 0 & A(\theta_4)  \\
 \end{array} \right) \in SO(8),
 \end{equation}

\noindent where
\begin{align}
& A(\theta_{\alpha})=
  \left( \begin{array}{cc}
  \cost_{\alpha} & \sit_{\alpha}  \\
  -\sit_{\alpha} & \cost_{\alpha}  \\
  \end{array} \right),  \qquad  \qquad
 \theta_{\alpha}=
p_{\alpha} t + q_{\alpha} s + l_{\alpha} r,
\end{align} 
with $t, s, r \in [ 0, 2\pi)$, and 
\begin{equation}
\Theta =  \left( \begin{array}{cccc}
p_1 & p_2&p_3&p_4  \\
q_1 & q_2&q_3&q_4  \\
l_1 & l_2&l_3&l_4  
  \end{array} \right) 
\end{equation}
is the matrix of the integral weights  defining the action.

Next, recall that the Hopf fibration $S^{31}  \longrightarrow \mathbb{H}P^7 $, acted on isometrically by $Sp(1)$ via left
 multiplication of quaternions,  gives as quotient:
\begin{equation}                  \label{C0 : 6}
 \left. \begin{array}{ccccc}
S^{31} &    \overset{Sp(1)}{\Longrightarrow} &
 \frac{SO(8)}{SO(4)\times Sp(1)}   &
 \overset{SO(3)}{\longrightarrow}  &  \frac{SO(8)}{SO(4)
 \times SO(4)}\cong \grm ,\\
 \end{array} \right.
\end{equation}
and we are going now to add to it the $T^3_{\Theta}-$ action.

Accordingly, we look at $G= Sp(1)\times T^3_{\Theta}$ as a
subgroup of the $3$-Sasakian isometries of $S^{31}$.
\noindent The moment maps $\mu : S^{31}
\rightarrow \mathfrak{sp}(1) \otimes \mathbb{R}^3\cong
\mathbb{R}^9$ associated with the $Sp(1)$ action
and: 
$ \nu : S^{31} \rightarrow \mathfrak{u}(1)^3\otimes \mathbb{R}^3$
associated with $T^3_{\Theta}\cong U(1)^3$
read respectively:
\begin{equation}           \label{C1: 9}
\mu (\uno)= ( \sum_{\alpha=1}^{8} \overline{u}_{\alpha} i u_{\alpha} \,\,
\sum_{\alpha=1}^{8} \overline{u}_{\alpha} j u_{\alpha}  \,\,
\sum_{\alpha=1}^{8} \overline{u}_{\alpha} k u_{\alpha} )\in \mathfrak{sp}(1) \otimes \mathbb{R}^3,
\end{equation}
and    

\begin{equation}             \label{C1: 10}
 \nu (\uno)= \left( \begin{array}{c}
 \sum_{\alpha = 1}^{4}  p_{\alpha}(\overline{u}_{2\alpha -1 } u_{2\alpha} -  \overline{u}_{2\alpha  } u_{2\alpha -1} )  \\
 \sum_{\alpha = 1}^{4}  q_{\alpha}(\overline{u}_{2\alpha -1 } u_{2\alpha} -  \overline{u}_{2\alpha  } u_{2\alpha -1} )  \\
 \sum_{\alpha = 1}^{4}  l_{\alpha}(\overline{u}_{2\alpha -1 } u_{2\alpha} -  \overline{u}_{2\alpha  } u_{2\alpha -1} )  \\
\end{array} \right)\in \mathfrak{u}(1)^3\otimes \mathbb{R}^3,
\end{equation}
 
\noindent where $\uno\in S^{31}\subset \mathbb{H}^8$.

The zero set $\mu^{-1}(0)$  can be easily
identified with the Stiefel manifold
of oriented orthonormal $4$-frames in $\mathbb{R}^8$, and it is therefore natural to 
look at elements of $N(\Theta)=\mu^{-1}(0)\cap\nu^{-1}(0)$ as $4\times 8$ real matrices
$ \uno =(u_1 , u_2 ,\dots, u_7 , u_8 )$,
whose columns $u_\rho$ are coefficients of a quaternion
respect to the base $\{
1,i,j,k\}$. Of course any such matrix $\boldsymbol{u}$
 has rank  $4$. 

\begin{defin} Let $\alpha = 1,2,3,4$. Any pair $\uua$  of 
quaternionic coordinates of
  $\uno\in S^{31}\subset \mathbb{H}^8$ will be called a \rm{quaternionic pair}.
\end{defin}

\begin{lem}
 Suppose that all the minor determinants
\begin{equation}  \label{C1: 10a}
\qquad \qquad \qquad \Delta_{\alpha\beta\gamma}=
\left\vert \begin{array}{ccc}
p_{\alpha} & q_{\alpha} & l_{\alpha} \\
p_{\beta} & q_{\beta} & l_{\beta} \\
p_{\gamma} & q_{\gamma} & l_{\gamma} \\
\end{array} \right\vert  \qquad \,\,\, \qquad \qquad
  (1\leq \alpha<\beta<\gamma\leq4) ,
\end{equation}
of $\Theta$ do not vanish. Then the zero set $N(\Theta)$  contains no elements $\uno$ having a null quaternionic pair.
\end{lem}
$Proof.$
Refer to the choice of  $(u_7 ,u_8 )$ as a null quaternionic pair on some point of $N(\Theta)$.
Let $x_{\alpha}= \overline{u}_{2\alpha -1 } u_{2\alpha} -  \overline{u}_{2\alpha  }
u_{2\alpha -1}, \,\, \alpha = 1,2,3,4$,
and rewrite $\nu (\uno )$ as 
\begin{equation}       \label{C1: 11}
\nu (\uno)= \Theta \left( \begin{array}{c}
x_1 \\
x_2 \\
x_3 \\
x_4 \\
\end{array}   \right)= \left( \begin{array}{ccc|c}
  &            &  & p_4\\
  & \mathbb{A} &  & q_4 \\
  &            &  & l_4 \\
\end{array} \right) \left( \begin{array}{c}
x_1 \\
x_2 \\
x_3 \\
x_4 \\
\end{array}   \right),
\end{equation}
where $\mathbb A={\left( \begin{array}{ccc}
p_1  & p_2 & p_3 \\
q_1  & q_2 & q_3 \\
l_1  & l_2 & l_3\\
\end{array} \right)}$.
Since $det\, \mathbb{A}= \Delta_{123} $, the equation $\nu (\uno )=0$
is solved by
$
\left( \begin{array}{c}
x_1 \\
x_2 \\
x_3 \\
\end{array} \right) =-\mathbb{A}^{-1}
 \left( \begin{array}{c}
 p_4 x_4 \\
  q_4 x_4\\
 l_4 x_4\\
\end{array} \right)
$, and
$N(\Theta)\cap \{ u_7= u_8 = 0\} $ has equations:  
\begin{equation} \label{C1: 14}
\left\{ \begin{array}{l}
\overline{u}_1 u_2 =\overline{u}_2 u_1 \\
\overline{u}_3 u_4 =\overline{u}_3 u_4 \\
\overline{u}_5 u_6 =\overline{u}_6 u_5 \\
\mu (\uno ) =\sum_{\alpha=1}^6 \overline{u}_{\alpha}\sigma u_{\alpha}|_{
\sigma = i,j,k}=0.
\end{array} \right.
\end{equation}
\noindent From (\ref{C1: 14}) we get 
$Im (\barun\um )=0,\,$  i.e. $\barun\um\in\mathbb{R}$.

Now observe that in equations (\ref{C1: 14}) we can assume, up to
a scale, that $u_{\alpha}$
belongs to $Sp(1)$. Thus 
maps
$
 (\un , \um ) \in Sp(1)\times Sp(1)\rightarrow \barun\um   \in Sp(1),$ $
\alpha = 1,2,3 $,
are consequently defined.
It follows
$\barun\um =\pm 1 $, $\un =\pm \um$,
and that  $\uno\in N(\Theta)$, as a a real $4\times 8 $  matrix, cannot
satisfy all the equations (\ref{C1: 14}):  the first three equations
force in fact the columns $\un={}^T(\un^0,\un^1,\un^2\un^3)$ and $\um=
{}^T(\um^0,\um^1,\um^2\um^3)$
of each quaternionic pair to be  proportional to each other.
 Thus the matrix
$\uno$ has at most rank $3$, contadicting the assumpion of $\uno$ as a $4$-frame
 in $\mathbb{R}^8$. 
It follows that  $N(\Theta)\cap \{ u_7= u_8 = 0\} $
is empty. $\square$

\begin{prop}
The action of $G=Sp(1)\times T^3_{\Theta}$ on $N(\Theta)$
 is locally free if and only if all
the following determinants:

\begin{align}  \label{C1: 15a}
& \qquad \qquad \qquad \qquad  \qquad \qquad \Delta_{\alpha\beta\gamma}=
\left\vert \begin{array}{ccc}
p_{\alpha} & q_{\alpha} & l_{\alpha} \\
p_{\beta} & q_{\beta} & l_{\beta} \\
p_{\gamma} & q_{\gamma} & l_{\gamma} \\
\end{array} \right\vert   \qquad \qquad
  (1\leq \alpha<\beta<\gamma\leq4) ,  
  \end{align}
  and:
  \begin{align}  \label{C1: 15b}
&  \qquad   \sq  :=
\left\vert \begin{array}{ccc}
p_1 \pm p_2 & q_1\pm q_2 & l_1 \pm l_2 \\
p_1 \pm p_3 & q_1\pm q_3 & l_1 \pm l_3 \\
p_1 \pm p_4 & q_1 \pm q_4 & l_1 \pm l_4 \\
\end{array} \right|   
\end{align}
do not vanish.
\end{prop}
$Proof$. By Lemma $2.1$  the conditions $\Delta_{\alpha\beta\gamma}\ne 0$ insure
that $N(\Theta)$ has no points with a null quaternionic pair.  Then the  fixed point
equations can be written as:  
\[
A(\theta_\alpha) \left( \begin{array}{c}
u_{2{\alpha}-1}  \\
u_{2{\alpha}}  \\
\end{array} \right) =
\left( \begin{array}{cc}
a_{\alpha} & b_{\alpha} \\
-b_{\alpha} & a_{\alpha} \\
\end{array} \right)  \left( \begin{array}{c}
u_{2{\alpha}-1}  \\
u_{2{\alpha}}  \\
\end{array} \right)  =
\lambda \left( \begin{array}{c}
u_{2{\alpha}-1}  \\
u_{2{\alpha}}  \\
\end{array} \right),
\]
where $a_{\alpha} = \cost_{\alpha}, \,\,
b_{\alpha} = \sit_{\alpha}$, and $\lambda\in Sp(1)$.
It follows:  
\begin{equation}      \label{C1: 18}
\left. \begin{array}{l}
a_{\alpha} |\un |^2 + b_{\alpha}\um \barun= \lambda |\un |^2 ,\\
-b_{\alpha} \un \barum+ a_{\alpha}  |\um |^2 = \lambda |\um |^2 ,\\
\end{array} \right.
\end{equation}
and:
\begin{equation}     \label{C1: 19}
\begin{aligned}
a_\alpha ( |\un |^2 + |\um |^2 ) +   b_\alpha (& \um \barun -  \un \barum ) = \lambda  (|\un |^2 + |\um |^2 ),
\end{aligned}
\end{equation} 
where by Lemma 2.1 the term multiplying $\lambda$ is non-zero. Also: 
\begin{equation}     \label{C1: 20}
Re  \; \lambda  =a_{\alpha},\,\,\,\,\,\, Im \; \lambda =
 b_{\alpha}\frac{ ( \um \barun -  \un \barum ) }{(|\un |^2 + |\um |^2 )},
\end{equation}
so that $   a_1 =a_2 = a_3= a_4 $
and $
b_1=\pm b_2 =\pm b_3=\pm b_4 .$
Therefore: 
 \begin{equation}         \label{C1: 23}
\left\{ \begin{array}{l}
(p_1 \pm p_2 )t + ( q_1 \pm q_2 )s + (l_1 \pm l_2 )r = 2h_{12}^{\pm} \pi  \\
(p_1 \pm p_3 )t + ( q_1 \pm q_3 )s + (l_1 \pm l_3 )r = 2h_{13}^{\pm}  \pi  \\
(p_1 \pm p_4 )t + ( q_1 \pm q_4 )s + (l_1 \pm l_4 )r = 2h_{23}^{\pm} \pi , \\
\end{array} \right.
\end{equation}
where $ h_{\alpha\beta}\in \mathbb{Z}$.
To have a locally free action, we need that all these eight systems have at most discrete solutions,
i.e. that
the eight determinants ${}^{1\pm 2}\square_{1\pm 4}^{1\pm 3}$     
do not vanish. $\square$

\begin{prop}
 There is no weight matrix $\Theta$ such that the action of
 $G=Sp(1)\times T^3_{\Theta}$ on $N(\Theta)$ is free.  
\end{prop}

$Proof.$ From the previous proof we see that there is a unique solution
 for the fixed point equations (\ref{C1: 23})
if and only if $\vert \sqpm \vert = 1.$  On the other hand, the identities:    

\begin{equation}     \label{C1: 26}
\left\{ \begin{array}{l}
{}^{1+ 2}\square_{1+ 4}^{1+ 3}= \Delta_{123}- \Delta_{124}+ \Delta_{134}+\Delta_{234}  \\
{}^{1+ 2}\square_{1+ 4}^{1- 3}= -\Delta_{123}- \Delta_{124}- \Delta_{134}-\Delta_{234} \\
{}^{1+ 2}\square_{1- 4}^{1+ 3}= \Delta_{123}+ \Delta_{124}-\Delta_{134}-\Delta_{234}  \\
{}^{1+ 2}\square_{1- 4}^{1- 3}= -\Delta_{123}+ \Delta_{124}+\Delta_{134}+\Delta_{234}  \\
{}^{1- 2}\square_{1+ 4}^{1+ 3}=- \Delta_{123}+ \Delta_{124}+ \Delta_{134}-\Delta_{234}  \\
{}^{1- 2}\square_{1+ 4}^{1- 3}= \Delta_{123}+ \Delta_{124}- \Delta_{134}+\Delta_{234}  \\
{}^{1- 2}\square_{1- 4}^{1+ 3}= -\Delta_{123}- \Delta_{124}- \Delta_{134}+\Delta_{234}  \\
{}^{1- 2}\square_{1- 4}^{1- 3}= \Delta_{123}- \Delta_{124}+ \Delta_{134}-\Delta_{234}  , \\
\end{array} \right.
\end{equation}
can be solved with respect to the $\Delta_{\alpha\beta\gamma}:$  

\begin{equation}     \label{C1: 27}
  {} \quad \left\{ \begin{array}{l}
 \Delta_{123} = -\frac{Y +W}{2}  \\
 \Delta_{124} = -\frac{X + Y}{2}\\
 \Delta_{134}= \frac{X +Y -Z + W}{2}\\
 \Delta_{234} =\frac{Z- Y}{2} ,
 \end{array}  \right.   
\end{equation}

\noindent where $X ={}^{1+ 2}\square_{1+ 4}^{1+ 3}$, $Y ={}^{1+ 2}\square_{1+ 4}^{1- 3}$, $ Z ={}^{1- 2}\square_{1- 4}^{1+ 3}$, $ W= {}^{1+ 2}\square_{1- 4}^{1- 3}$.
In particular:   
\begin{equation}          \label{C1: 29}
\left\{ \begin{array}{l}
\pm 1 = {}^{1+ 2}\square_{1- 4}^{1+ 3}=  -(X + Y + W) \\
\pm 1 = {}^{1- 2}\square_{1+ 4}^{1- 3}=  Z - X -2Y -W \\
\pm 1 = {}^{1- 2}\square_{1+ 4}^{1+ 3}= Y  + Z + W\\
\pm 1 = {}^{1- 2}\square_{1- 4}^{1- 3}= X + Y - Z, \\
\end{array}\right.
\end{equation}

\noindent  and our assumptions $\Delta_{\alpha\beta\gamma} \neq 0$ and $X,Y,Z,W=\pm 1$  give $
X = Y =W =-Z$, so that $(X, Y, Z, W )$ is either $(1,1,-1,1)$ or $(-1, -1, 1, -1)$
and $(\Delta_{123} , \Delta_{124}, \Delta_{134}, $
$\Delta_{234})$ is either
 $(1,1, -2, 1)$ or $(-1, -1, 2, -1)$.
Since these choices of $(X,Y,Z,W)$
are not a solution
for the above system, no free action can be obtained. $\square$ \\

We conclude the paragraph by summarising all of this in a statement. Note that
Theorem A of the Introduction then follows.

\begin{teor}
The action of $Sp(1)\times T^3_{\Theta}$ on $N(\Theta)=\nu^{-1}(0)\cap
 \mu^{-1}(0)$ is never free, and it is locally free if and only if the following conditions hold:     
\bigskip

i) $\Delta_{\alpha\beta\gamma}\ne 0$\,\,\,\, for any
$(\alpha,\beta,\gamma)$,
\bigskip

ii) all the determinants ${}^{1\pm 2}\square_{1\pm 3}^{1\pm 4} $
  are non zero.
\bigskip

In such a case the quotient
\begin{equation}       \label{C1: 29a}
\mathcal{M}^7(\Theta)= \frac{N(\Theta)}{Sp(1)\times T^3_{\Theta}}
\end{equation}
 is a compact $7-$dimensional $3-$Sasakian orbifold and a principal $SO(3)$-bundle over a $4$-dimensional orbifold $\mathcal{O}^4(\Theta )$ with a positive SDE metric
and a one-dimensional group of isometries.
\end{teor}

The orbifolds $\mathcal{M}^7(\Theta)$ are not toric. To see this,
look at the foliation on 
 $N(\Theta)$ that gives any such orbifold as 3-Sasakian quotient. Then observe that $N(\Theta)$ is a compact submanifold of  $S^{31}\subset
  \mathbb{H}^{8}$ as the zero locus of the
quadratic functions defined by the moment maps $\mu$ and $\nu$.
Thus all the isometries
of $N(\Theta)$ come from the restriction of the isometries of $S^{31}$ and, projecting to the $4$-dimensional base,
the group of isometries associated
to $\mathcal{O}^4(\Theta )$ turns out to be one-dimensional. 

\vspace{0.3 cm}

\begin{flushleft}
 $\mathbf{Examples\,\, 2.2.}$
There are many matrices which satisfy the assumptions of  Theorem $2.1$ and hence of Theorem A in the Introduction. For example:
 \begin{equation}    \label{C1: 30}
\Theta_1 = \left( \begin{array}{cccc}
1 & 0 & 1 & 1 \\
0 & 1 & 1 & 1  \\
1 & 1 & 0 & 1  \\
\end{array} \right), \qquad \Theta_2 = \left( \begin{array}{cccc}
9 & 2 & 7 & 1 \\
40 & 9 & 31 & 0  \\
1 & 2 & 0 & 1  \\
\end{array} \right)
\end{equation}

\end{flushleft}

\bfl have minor determinants $\Delta_{\alpha \beta \gamma} = (-2, -1,1,-1)$ and $\Delta_{\alpha \beta \gamma} = (1,72,-32,-63)$, respectively. All conditions listed in Theorem $2.1$ are easily verified. $\square$ \efl

\vspace{1 cm}

\section{The Singular Locus of $\mathcal{Z}^6(\Theta)$ }

Let $\mathcal{Z}^6(\Theta)$ be the twistor space of any of the orbifolds $\mathcal{O}^4(\Theta )$ constructed in Section 2 and let $\Sigma(\Theta)$ be its singular locus .  The zero set $N(\Theta)=   $
 $\mu^{-1}(0) \cap \nu^{-1}(0) \subset S^{31}$ is acted on by the group $\widetilde{G}=G\times U(1)=T^3_{\Theta}\times Sp(1)\times U(1)$ that, up to the central $\mathbb{Z}_2$, is
 a subgroup  of $Sp(8)\cdot Sp(1)\subset SO(32)$. Let:
 \begin{equation}    \label{C2 : 1}
  \Phi:\, T^3_{\Theta}\times Sp(1)\times U(1)\times
  N(\Theta) \longrightarrow N(\Theta)
  \end{equation}
be the action, where: 
 \begin{equation} 
\Phi \big( ( A(\Theta), \lambda, \rho)
 \big)\big((\underline{z}, \underline{w} ) \big)
 =   A(\Theta) \lambda  \left( \begin{array}{c}
\underline{z}\\
\underline{w} \\
\end{array}\right)
 \rho,
 \end{equation}
and we have identified $\mathbb{H}^8 \cong \mathbb{C}^8\times \mathbb{C}^8$ by
  $u_{\alpha} = z_{\alpha} + jw_{\alpha}$. Thus
  $\zw =\boldsymbol{u}=
 (u_1,u_2,\ldots,u_7,u_8)\in \mathbb{H}^8$ and we will
  use both notations $\uno$ and $\zw$.
The twistor space $\mathcal{Z}^6(\Theta)$ 
is the leaf
space of the $\widetilde{G}-$action on $N(\Theta)$. There is a natural
stratification of $\mathcal{Z}^6(\Theta)$
 and we want to see how any singular stratum in
   $\mathcal{Z}^6(\Theta)$
appears from the action of $\widetilde{G}$ on $N(\Theta)$.

\begin{defin}
 We say that two points $\zw$, $(\underline{z}_1, \underline{w}_1)\in
 N(\Theta)$ define the same $\widetilde{G}-$stratum $\overline{S}$ of $\mathcal Z^6(\Theta)$ if their corresponding
 isotropy subgroups $\widetilde{G}_{\zw}$, $\widetilde{G}_{(\underline{z}_1,
 \underline{w}_1)}$   are conjugate with respect to the $\widetilde{G}$-action.
\end{defin}

To get the possible isotropy subgroups, fix a point
$\zw=\uno$ and write 
 the fixed point equations in quaternionic pairs $\uu$, as follows:

  \begin{equation} \label{C2a : 1}
 A(\theta_{\alpha})
\left( \begin{array}{c}
u_{2\alpha-1} \\
u_{2\alpha}  \\
 \end{array}   \right) =
\left( \begin{array}{c}
\lambda u_{2\alpha-1} \rho \\
\lambda u_{2\alpha} \rho \\
 \end{array}   \right),
\end{equation}
where $A(\theta_{\alpha})=\left(\begin{array}{cc}
\cost_{\alpha} & \sit_{\alpha} \\
-\sit_{\alpha}  & \cost_{\alpha} \\
\end{array} \right)$,
$\lambda= \epsilon + j\sigma\in Sp(1)$, $\rho\in U(1)$.
Equivalently: 
\begin{equation}  \label{C2 : 8}
A(\theta_{\alpha})
 \left(\begin{array}{cc}
z_{2{\alpha}-1} & w_{2{\alpha}-1} \\
z_{2{\alpha}} & w_{2{\alpha}} \\
\end{array} \right) =
{\left[ \begin{array}{c}
\left(\begin{array}{cc}
\epsilon & -\overline{\sigma}  \\
\sigma & \overline{\epsilon} \\
\end{array} \right)
 \left(\begin{array}{cc}
z_{2{\alpha}-1} & z_{2{\alpha}} \\
w_{2{\alpha}-1} & w_{2{\alpha}} \\
\end{array} \right)
 \left(\begin{array}{cc}
\rho& 0 \\
0 & \rho \\
\end{array} \right)
\end{array}\right]}^T .
\end{equation}

\noindent The following property is easily verified:

\begin{lem}
Let $\zw$ be a point in $N(\Theta)$.
 Then, up to $\widetilde{G}-$conjugation, we have
 $\widetilde{G}_{\zw}\subset T^3_{\Theta}\times U(1)^{\epsilon} \times U(1)$, where $ U(1)^{\epsilon}= \{\lambda\in Sp(1)\,\, |\,\,
 \sigma= 0 \}$.
\end{lem}

Thus, orbits through points
  $\zw\in N(\Theta)$
 with non trivial isotropy subgroup
   $\widetilde{G}_{\zw}$
   $\subset T^3_{\Theta}\times U(1)^{\epsilon}\times U(1)$ give all the $\widetilde{G}-$strata of $N(\Theta)$ whose projection gives rise to singular strata of $\mathcal Z^6(\Theta)$.

Rewrite now equations (\ref{C2 : 8}) as follows ($\alpha=1,2,3,4$):
\begin{equation}    \label{C2 : 9}
\overbrace{\left( \begin{array}{cc|cc}
0                          & -\sigma\rho &  -\sit_{\alpha}  & \cost_{\alpha} - \overline{\epsilon}\rho  \\
-\sigma\rho     & 0                      & \cost_{\alpha} - \overline{\epsilon}\rho & \sit_{\alpha} \\
\hline
-\sit_{\alpha}   & \cost_{\alpha} - \epsilon\rho           &    0       &      \overline{\sigma}\rho  \\
\cost_{\alpha} -\epsilon \rho & \sit_{\alpha} &  \overline{\sigma}\rho   &   0    \\
\end{array} \right)}^{M_{\alpha}:=}
\left(\begin{array}{c}
z_{2\alpha-1} \\
z_{2\alpha}\\
w_{2\alpha-1} \\
w_{2\alpha} \\
\end{array} \right) =
\left(\begin{array}{c}
0 \\
0\\
0\\
0 \\
\end{array} \right),
\end{equation}
\noindent and note that
none of the $M_{\alpha}$
can have  rank $4$, since
 otherwise the correspondent quaternionic
  pair $(u_{2\alpha-1}, u_{2\alpha})$ would vanish, 
a contradiction with Lemma $2.1$.

\begin{prop}
Let $M_{\alpha}$ be the matrix  in formula $(\ref{C2 : 9})$. Then
 $det\, \mal=0$ if and only if at least one of the following four identities
\begin{equation}   \label{C2 : 18b}
\overline{\rho}e^{\pm i\theta_{\alpha}}= Re \; \epsilon \pm  i\sqrt{(Im \: \epsilon)^2 + |\sigma|^2}
\end{equation}
holds.
\end{prop}

$Proof.$
By using the block notation: 
\begin{equation}   \label{C2 : 10}
\begin{aligned}
M_{\alpha}&=\left(\begin{array}{c|c}
A & B \\
\hline
 C & D \\
 \end{array}\right)
 \end{aligned}
\end{equation}
one has, for $\sigma \neq 0$:
\begin{equation}   \label{C2 : 11}
M_{\alpha}=  \left(\begin{array}{cc}
 B & Id\\
  D & 0 \\
 \end{array}\right) \left(\begin{array}{cc}
 D^{-1}C &  Id \\
 A- BD^{-1}C  &  0\\
 \end{array}\right),
\end{equation}
where the matrix $ (A- BD^{-1}C)$ is given by:
\[   
\frac{\sigma}{|\sigma|^2}\left(\begin{array}{cc}
\overline{\rho} sin\, 2\theta_{\alpha} -2 \sit_{\alpha} Re \; \epsilon &   ( 2 \cost_{\alpha} Re\; \epsilon- \overline{\rho}cos\, 2\theta_{\alpha}) - \rho \\
( 2 \cost_{\alpha} Re \; \epsilon- \overline{\rho}cos\, 2\theta_{\alpha}) - \rho  &   -\overline{\rho} sin\, 2\theta_{\alpha} +2 \sit_{\alpha} Re \; \epsilon \\
 \end{array} \right).
\]
It follows:
 \begin{equation} \label{C2 : 16}
\begin{aligned}
 &det\, \mal = \overline{\sigma}^2\rho^2 det\, (A- BD^{-1}C )= \\
 &=\rho^2\bigg( \rho+
 \overline{\rho} (e^{-i\theta_{\alpha}})^2
- 2 (Re \;\epsilon) e^{-i\theta_{\alpha}}\bigg)\bigg( \rho+
\overline{\rho} (e^{i\theta_{\alpha}})^2
- 2 (Re \;\epsilon) e^{i\theta_{\alpha}}\bigg),
\end{aligned}
\end{equation}
 and this is zero if and only if: 
\begin{equation}  \label{C2 : 17}
  \overline{\rho} (e^{-i\theta_{\alpha}})^2 - 2 (Re \;\epsilon)
   e^{-i\theta_{\alpha}}+ \rho=0
  \,\,\,\,\,\, \mathrm{or}\,\,\,\,\,\,
  \overline{\rho} (e^{i\theta_{\alpha}})^2- 2 (Re\;\epsilon)
  e^{i\theta_{\alpha}}+ \rho=0,
\end{equation}
that gives the stated condition. 
 $\square$ 
 \newline
\begin{oss} Note that, when $\sigma = 0$, the condition $det\, M_{\alpha}=0$ is
 equivalent to
\begin{equation}  \label{C2 : 16u}
\epsilon\rho= e^{\pm i\theta_{\alpha}}\,\,\,\,  or \,\,\,\,  \overline{\epsilon}\rho=
 e^{\pm i\theta_{\alpha}},
\end{equation}
which are special cases of formula (\ref{C2 : 18b}).
\end{oss}

We can rephrase all of this as follows:

\begin{prop}
\bfl For any $(\underline{z}, \underline{w})\in N(\Theta)\subset S^{31}$, to get a non trivial solution for the fixed
point equations $(\ref{C2a : 1})$ it is necessary that condition $(\ref{C2 : 18b})$ holds for some choices of the signs and for $\alpha=1,...,4$. \efl
\end{prop}

Assume now $\sigma = 0$ in system $(\ref{C2 : 9})$ and use in each block $M_{\alpha}$
 one or two relations among the four in (\ref{C2 : 16u}). Then we see that equations
 (\ref{C2 : 9}) admit non null solutions $\zw\in
 \mathbb{H}^8$,  fixed by a subgroup
 $H_{((\rho\epsilon, \rho\overline{\epsilon}), \alpha)}$ of the group generated by the chosen
 relations. Thus, Proposition $3.2$ gives that for any of these solutions
 $\zw$
 the isotropy subgroup $\widetilde{G}_{\zw}$
is contained in  $H_{((\rho\epsilon, \rho\overline{\epsilon}), \alpha)}$.
Depending on the numbers of the relations (\ref{C2 : 16u}),
the following possibilities for
 the rank of the blocks $M_{\alpha}$ can occur:

\vspace{0,5 cm} 
\bfl
$1)$\,\,\,  just one of the relations in  (\ref{C2 : 16u}) holds $\iff$
$rank\, \mal=3,$ \\
\vspace{0,2 cm}
$2)$\,\,\,   two relations in (\ref{C2 : 16u}) hold $\iff$  $rank\,
\mal=2,$  \\
\vspace{0,2cm}
$3)$\,\,\,  three or four relations in (\ref{C2 : 16u}) are satisfied
$\iff$  $ \mal=0_{4\times 4}.$
\efl
\vspace{0,5 cm}
When three or four
 relations  in  (\ref{C2 : 18b})
 hold, for each $\alpha=1,2,3,4$, they describe
 the non effectivity. Thus, up to the non effective subgroup,
the third case can be ignored. Accordingly:

\begin{lem}
Assume $\sigma = 0$ in system $(\ref{C2 : 9})$. If $rank\,M_{\alpha}= 3$, 
its solutions are given by any of the following:
 \begin{equation}  \label{C2 : 30}
\begin{aligned}
\qquad \qquad &1)\,\,\,   {}^{\pm}V_1^{\alpha} =\{(z_{2\alpha-1}, \pm iz_{2\alpha-1}, 0,0)\}, \quad\quad  \rho \epsilon = e^{\pm i \theta_\alpha},\\
&    \\
\qquad \qquad & 2)\,\,\, {}^{\pm}V_2^{\alpha} =\{( 0,0, w_{2\alpha-1}, \pm iw_{2\alpha-1})\}, \quad\quad  \rho \bar \epsilon = e^{\pm i \theta_\alpha},   {}
\end{aligned}    \hspace{2 cm}
\end{equation}  \\
and when $rank\,M_\alpha =2 $ by any of :
\begin{equation}   \label{C2 : 30b}
\begin{aligned}
\qquad  & 3)\,\,\, {}^{(\pm, \pm )}V_3^{\alpha} =\{(z_{2\alpha-1}, \pm iz_{2\alpha-1}, w_{2\alpha-1}, \pm iw_{2\alpha-1} )\}, \quad  \rho \epsilon = e^{\pm i \theta_\alpha}, \rho \bar \epsilon = e^{\pm i \theta_\alpha}, \\
& \\
\qquad  & 4)\,\,\,  V_4^{\alpha} =\{(z_{2\alpha-1}, z_{2\alpha }, 0,0)\},  \rho  \epsilon = e^{ i \theta_\alpha} = e^{- i \theta_\alpha}, \\
& \\
\qquad &5)\,\,\,V_5^{\alpha} =\{( 0,0,
w_{2\alpha-1},
w_{2\alpha})\}. \rho \bar \epsilon = e^{i \theta_\alpha}  = e^{-i \theta_\alpha}.{}
 \end{aligned}
 \end{equation}
\end{lem}

\vspace{0,3 cm}

\noindent It follows:

\begin{cor}
 Let $\zw$ be a point in $N(\Theta)$ with non trivial isotropy
subgroup $\widetilde{G}_{\zw}$.
Then each quaternionic pair $\uu$ of $\zw$
belongs to one of the sets $\widetilde{G}\cdot {}^{\pm}V_{1}^{\alpha}$,
$\widetilde{G}\cdot {}^{\pm}V_{2}^{\alpha}$, $\widetilde{G}\cdot {}^{(\pm, \pm)}V_{3}^{\alpha}$,
$\widetilde{G}\cdot V_{4}^{\alpha}$ or $\widetilde{G}\cdot V_{5}^{\alpha}$. 
\end{cor}

\begin{prop}
Let $\zw$ be  a point on a singular 
$\widetilde{G}-$orbit of $N(\Theta)$, and assume that the hypotheses of Theorem A hold.
Then at least one of the blocks 
$ M_{\alpha}$,  $\alpha=1,2,3,4$, has $rank=2$.  
\end{prop}

 $Proof.$
Assume that
$rank\, M_{\alpha}=3$ for
$\alpha=1,2,3,4$. 
Then, for each
$M_{\alpha}$,
just one of relations 
(\ref{C2 : 16u}) holds. Note first that there are  $\gamma \neq \delta$
such that $M_\gamma$ satisfies one of the first two identities in (\ref{C2 : 16u}) and $M_\delta$ one of the remaining two. In particular,
$e^{i\theta_{\gamma}} \neq
 e^{\pm i \theta_{\delta}}$ since otherwise
$M_{\gamma}$ and $M_{\delta}$ would have
 $rank\,= 2$.  In fact,
assuming that such indices $\gamma$ and $\delta$ do not exist, 
then solutions for
equations (\ref{C2 : 9}) would
have quaternionic pairs
either contained in a
${}^{\pm}V^{\alpha}_1$
or in a ${}^{\pm}V^{\alpha}_2$ for all $\alpha$.
Thus, these
solutions 
would be 
$8\times 4$ real matrices with $rank < 4$,
and as such
not points in $ N(\Theta)$.
Thus, let $M_{\gamma}$ and $M_{\delta}$ with the mentioned property, so with spaces of solutions 
of type ${}^{\pm}V^{\gamma}_1$ and
${}^{\pm}V^{\delta}_2$. Then, looking at all the indices $\alpha$ we get solutions in any of the following subspaces:
 \begin{equation}  \label{C2 : 34}
  \begin{aligned}
& \Sigma^{123}_4 =\Bigg\{ X\in M_{2\times8}
(\mathbb{C})\,\, |\,\, X=\left(\begin{array}{cc|cc|cc|cc}
* & *  & * & * & * & * & 0 & 0\\
0 & 0  & 0 & 0 & 0 & 0 & * & * \\
\end{array}\right)\Bigg\}, 
  \end{aligned}
 \end{equation}
  \begin{equation}  \label{C2 : 34a}
 \begin{aligned}
& \Sigma{}^{12}_{34} =\Bigg\{ X\in M_{2\times8}
(\mathbb{C})\,\, |\,\, X=\left(\begin{array}{cc|cc|cc|cc}
* & *  & * & * & 0 & 0 & 0 & 0\\
0 & 0  & 0 & 0 & * & * & * & *\\
\end{array}\right)\Bigg\}, 
 \end{aligned}
 \end{equation}
or in the similarly defined $ \Sigma^{124}_3,
 \Sigma^{134}_2,  \Sigma^{234}_1,
 \Sigma^{1}_{234},
 \Sigma^2_{134},  \Sigma^3_{124},\Sigma^4_{123}$, 
or in $\Sigma^{13}_{24},
\Sigma^{14}_{23}, \Sigma^{23}_{14}, 
\Sigma^{24}_{13}, \Sigma^{34}_{12}.$
In these matrices, the $2\times 2$
blocks  represent elements either
of ${}^{\pm}V^{\alpha}_1$ or of ${}^{\pm}V^{\alpha}_2$. 
Look for example at $ \Sigma^{123}_4$
with the choices (following notations in \ref{C2 : 30}):
$\uu\in{}^{+}V_1^{\alpha}$,  
$\alpha=1,3$, $(u_3,u_4)\in{}^{-}V_1^{2}$
and $(u_7,u_8)\in{}^{+}V_2^{4}$.
 By reading the $T^3_{\Theta}-$moment
 map equations on this set of solutions, we get:
 \begin{equation}       \label{C2a : 72a}
 \left\{\begin{array}{l}
p_1|z_1|^2- p_2|z_3|^2 +p_3|z_5|^2 +p_4|w_7|^2=0\\
 q_1|z_1|^2- q_2|z_3|^2 +q_3|z_5|^2 +q_4|w_7|^2=0\\
 l_1|z_1|^2- l_2|z_3|^2 +l_3|z_5|^2 +l_4|w_7|^2=0.\\
 \end{array}\right. \hspace{ 3 cm} {}
 \end{equation}
Similarly, the $Sp(1)-$moment map equation for this choice gives:
 \begin{equation}     \label{C2a : 73a}
 |z_1|^2  +|z_3|^2 +|z_5|^2 -|w_7|^2=0,
  \end{equation}
and we can rewrite all these 
equations as follows:  
 \begin{equation}      \label{C2a : 74a}
\left( \begin{array}{cccr}
  p_1 & -p_2 & p_3  & p_4 \\
  q_1 & -q_2 & q_3  & q_4 \\
  l_1 & -l_2 & l_3  & l_4 \\
  1 & 1  & 1 & -1 \\
 \end{array} \right)
 \left( \begin{array}{c}
  |z_1|^2 \\
  |z_3|^2 \\
  |z_5|^2 \\
  |w_7|^2 \\
 \end{array} \right)    =
 \left( \begin{array}{c}
  0\\
  0\\
  0\\
  0\\
 \end{array} \right),
 \end{equation}
where the determinant is
one of  the $\sq\ne 0$. Thus we
 get only the trivial solution.
  All the other listed cases can be treated similarly. $\square$

\begin{prop}
 Let $(\underline{z}, \underline{w})$ be
on a singular $\widetilde{G}-$orbit of $N(\Theta)$ and
assume that $M_{\beta}$ has $rank  =2$ for some $\beta =1,...,4$  with solutions of the corresponding system of type
 $3)$ in (\ref{C2 : 30b}).
 Then all the four matrices $M_{\alpha}$ have $rank \leq 2$.   
\end{prop}

 $Proof.$
 Without loss of generality, we can assume $\beta=1$ and that the
 conditions in (\ref{C2 : 16u}) giving solutions of type
 $3)$ in (\ref{C2 : 30b}) are:
 \begin{equation}    \label{C2 : 38a}
 \left\{\begin{array}{l}
\epsilon \rho= e^{i\theta_1}\\
\overline{\epsilon} \rho= e^{i\theta_1}.\\
\end{array} \right.
 \end{equation}
 Since
$rank\,\ M_{\alpha}\leq 3$ ($\alpha > 1$), at least one of relations (\ref{C2 : 16u}) holds, so that
either $\epsilon\rho = e^{\pm i\theta_{\alpha}}$ or
 $\overline{\epsilon}\rho = e^{\pm i\theta_{\alpha}}$. If one of the first two identities holds we get either  
$
\epsilon\rho = \overline{\epsilon}\rho =  e^{i\theta_1} = e^{ i\theta_{\alpha}}$ or
$
\epsilon\rho = \overline{\epsilon}\rho =  e^{i\theta_1} = e^{ -i\theta_{\alpha}}$, and thus
we also get one of the further conditions $\overline{\epsilon}
 \rho =e^{\pm i\theta_{\alpha}} $, so that $rank\, M_{\alpha} \leq 2$ for
 $\alpha=2,3,4$. $\square$

\vspace{0,5 cm}
It follows:

\begin{cor}
Assume the hypotheses of Theorem A, and let $\zw\in N(\Theta)$ be such that one of its quaternionic pairs
 belongs to a ${}^{(\pm, \pm )}V_3^\alpha$,
and another one 
either to a $V^\beta_4$ or to a $V^\beta_5$.
Then the isotropy subgroup $\widetilde{G}_{\zw}$ is trivial.  Then,
 if $\zw\in N(\Theta)$ is a singular point with  a quaternionic
 pair $\uu$ contained in a ${}^{(\pm,\pm)}V_3^\alpha$,
then all its quaternionic pairs are contained in one of the spaces
${}^{(\pm ,\pm )}V_3^\beta$.    
\end{cor}

\vspace{0,3 cm}
 
Lemma $3.2$ can be viewed as a description of 
strata on $\mathcal Z^6(\theta)$ coming from the action of $\widetilde{G}$ on $N(\Theta)$.
In fact, quaternionic pairs corresponding to singular points in the quotient are listed in (\ref{C2 : 30})
and (\ref{C2 : 30b}). In particular, by Propositions $3.3, 3.4$ and Corollary $3.2$ we see that singular strata on $\mathcal Z^6(\Theta)$ can be distinguished into the following two different families.
The first family of singular strata on $\mathcal Z^6(\Theta)$ comes from 
points $\zw$ such that all of their quaternionic pairs are contained
in ${}^{(\pm ,\pm )}V_3^\alpha$, $\alpha=1,2,3,4$.  In the second family,
the $\zw$ have no quaternionic pairs contained in ${}^{(\pm ,\pm )}V_3^\alpha$.

\vspace{0,5 cm}

We begin now by studying  the first mentioned family of singularities on $\mathcal Z^6(\Theta)$. Here all quaternionic pairs of $\zw$ are in a ${}^{(\pm ,\pm )}V_3^\alpha$, so that:
\begin{equation}   \label{C2 : 41}
{}^T(\underline{z}, \underline{w} ):=
\left(\begin{array}{cc|cc|cc|cc}
z_1 & \pm iz_1  & z_3 & \pm iz_3 & z_5 & \pm iz_5 & z_7 & \pm iz_7\\
w_1 & \pm iw_1 &  w_3 & \pm iw_3 & w_5 & \pm iw_5 & w_7 & \pm iw_7 \\
\end{array}\right),
\end{equation}
where all the signs can be chosen independently.
The fixed point equations are:
\begin{equation}       \label{C2 : 42}
\left\{ \begin{array}{l}
e^{i(\theta_1\pm\theta_\alpha)}=1\\
   \epsilon\rho = e^{\pm i\theta_1}\\
\overline{\epsilon}\rho = e^{\pm i\theta_1}, \\
\end{array} \right.
\end{equation}
where $\alpha=1,2,3,4$. It is convenient to introduce the following notation.  
Fix a pair of
 signs $(\pm, \pm)$  in relations
  $\left\{\begin{array}{l}
\epsilon\rho = e^{\pm i\theta_1}, \\
\overline{\epsilon}\rho = e^{\pm i\theta_1}, \\
\end{array}\right.$
 and a  triple of signs
$(\pm, \pm, \pm)$,
in  $e^{i\theta_1}= e^{\pm  i\theta_{\alpha}}$, $\alpha=2,3,4$.
Any space of solutions of (\ref{C2 : 9}) is associated to a  $5-$tuple of signs $\boldsymbol{\pm} = \big( (\pm,\pm ),(\pm,\pm,\pm ) \big)$.  Accordingly, we will denote any such space  of solutions by
  ${S}{}^{\boldsymbol{ \pm}}$.

Next, consider the intersections of all spaces 
  ${S}{}^{\boldsymbol{ \pm}}$ with $N(\Theta)$. A first observation is the following.

\begin{prop}  Let $\boldsymbol{ \pm} = \big( (+,- ),(\pm,\pm,\pm ) \big)$ or $\boldsymbol{ \pm} = \big( (-,+),(\pm,\pm,\pm ) \big)$. Then
  ${S}{}^{\boldsymbol{ \pm}}$ has an empty intersection with
$N(\Theta)=\mu^{-1}(0)\cap \nu^{-1}(0)$.
\end{prop}

 $Proof.$  It is sufficient to look at the intersection with 
 ${S}{}^{(+,-)}_{(\pm, \pm, \pm)}$ (the other case
is symmetric). By reading the $Sp(1)$-moment map equations on points $\zw$ in
 ${S}{}^{(+,-)}_{(\pm, \pm, \pm)}$, we see that
 
\begin{equation}      \label{C2 : 49}
\left\{\begin{array}{l}
\sum_{\alpha=1}^4  (|z_{2{\alpha}-1}|^2 -|w_{2\alpha-1}|^2)=0 ,\\
\sum_{\alpha=1}^4  z_{2\alpha-1}w_{2\alpha-1}=0. \\
\end{array}\right.
\end{equation}

\noindent The moment map $\nu$ of the $T^3_{\Theta}$-action yields:

\begin{equation}            \label{C2 : 50}
\left\{\begin{array}{l}
\sum_{\alpha=1}^4 d_{\alpha} Im(z_{2\alpha-1}\overline{z}_{2\alpha} +
w_{2\alpha-1}\overline{w}_{2\alpha})=0, \\
\sum_{\alpha=1}^4 d_{\alpha} (z_{2\alpha-1}w_{2\alpha}-  z_{2\alpha} w_{2\alpha-1} ) =0,
\end{array}\right.
\end{equation}

\noindent  where $d_{\alpha} = p_{\alpha}, q_{\alpha}, l_{\alpha}$,
$\alpha=1,2,3,4$. Thus for points in ${S}{}^{(+,-)}_{(\pm, \pm, \pm)}$ :

\begin{equation}         \label{C2 : 51}
\left\{\begin{array}{l}
\sum_{\alpha=1}^4  (-1)^{m_{\alpha}} d_{\alpha}(|z_{2\alpha-1}|^2 -|w_{2\alpha-1}|^2)=0 ,\\
\sum_{\alpha=1}^4 (-1)^{m_{\alpha}} d_{\alpha} z_{2\alpha-1}w_{2\alpha-1}=0, \\
\end{array}\right.
\end{equation}

\noindent where the indices $m_{\alpha}$ depend on the $5$-ples of signs. If $\Gamma_{\alpha}:=|z_{2\alpha-1}|^2 -|w_{2\alpha-1}|^2$, we can rewrite all our equations in (\ref{C2 : 49}) and (\ref{C2 : 50}) as

\begin{equation}         \label{C2 : 52}
 \left\{\begin{array}{l}
\sum_{\alpha=1}^4 (-1)^{m_{\alpha}} p_{\alpha} \Gamma_{\alpha}=0 \\
\sum_{\alpha=1}^4 (-1)^{m_{\alpha}} q_{\alpha} \Gamma_{\alpha}=0 \\
\sum_{\alpha=1}^4 (-1)^{m_{\alpha}} l_{\alpha} \Gamma_{\alpha}=0 \\
 \sum_{\alpha=1}^4  (\Gamma_{\alpha})=0 \\
\end{array}\right.\,\, , \,\,
\left\{\begin{array}{l}
 \sum_{\alpha=1}^4 (-1)^{m_{\alpha}} p_{\alpha} z_{2\alpha-1}w_{2\alpha-1}=0\\
\sum_{\alpha=1}^4 (-1)^{m_{\alpha}} q_{\alpha} z_{2\alpha-1}w_{2\alpha-1}=0\\
\sum_{\alpha=1}^4 (-1)^{m_{\alpha}}l_{\alpha} z_{2\alpha-1}w_{2\alpha-1}=0 \\
\sum_{\alpha=1}^4 z_{2\alpha-1}w_{2\alpha-1}=0
\end{array}\right.
\end{equation}
and we can observe that the first four equation have the same determinant of coefficients as the last four, namely one of  the
 $\sq\ne0$.  Then:
 \begin{equation}     \label{C2 : 53}
 |z_{2\alpha-1}|^2 =|w_{2\alpha-1}|^2>0\,\,\, {\rm and} \,\,\,\, z_{2\alpha-1}w_{2\alpha-1}=0,
 \end{equation}
so that our system does not admit solutions. $\square$

\vspace{0.3 cm}

The following Proposition shows the existence in the singular locus of the two $2$-spheres appearing in Theorem B (i) of the Introduction.

\begin{prop}
Just one among the spaces
 ${S}{}^{(+,+)}_{(\pm, \pm, \pm)}$ and just one among the
  ${S}{}^{(-,-)}_{(\pm, \pm, \pm)}$
intersects 
   $N(\Theta)$.
\end{prop}

 $Proof.$ We outline the argument for the first set of spaces ${S}{}^{(+,+)}_{(\pm, \pm, \pm)}$, the other case being very similar.
 On any of the  ${S}{}^{(+,+)}_{(\pm, \pm, \pm)}$ the $Sp(1)-$
 moment map equation is given by:
 \begin{equation}   \label{C2 : 55}
  \left\{\begin{array}{l}
 \sum_{\alpha=1}^4  (|z_{2\alpha-1}|^2 -|w_{2\alpha-1}|^2)=0 ,\\
 \sum_{\alpha=1}^4  \overline{w}_{2\alpha-1}z_{2\alpha-1}=0, \\
 \end{array}\right.
 \end{equation}
representing the Stiefel manifold $\mathcal{S}=U(4)/U(2)$.
The $T^3_{\Theta} $ moment map equation depends instead on the chosen  ${S}{}^{(+,+)}_{(\pm, \pm, \pm)}$:
 \begin{equation}     \label{C2 : 56}
\sum_{\alpha=1}^4  (-1)^{m_{\alpha}}d_{\alpha}(|z_{2\alpha-1}|^2
 +|w_{2\alpha-1}|^2)=0,
 \end{equation}
with $d_{\alpha} =p_{\alpha}, q_{\alpha}, 
l_{\alpha}$.
This can be rewritten in quaternionic coordinates: 
 \begin{equation}    \label{C2 : 57}
 \sum_{\alpha=1}^4  (-1)^{m_{\alpha}}d_{\alpha} 
|u_{2\alpha-1}|^2=0,
\end{equation}
where $u_{2\alpha-1} =z_{2\alpha-1}+jw_{2\alpha-1}$, and
by intersecting with the 
 sphere $S^{31}$ we get:
\begin{equation}         \label{C2 : 59}
\left(\begin{array}{cccc}
(-1)^{m_1}p_1  & (-1)^{m_2}p_2 & (-1)^{m_3} p_3 & (-1)^{m_4} p_4 \\
(-1)^{m_1}q_1  & (-1)^{m_2}q_2 & (-1)^{m_3} q_3 & (-1)^{m_4} q_4 \\
(-1)^{m_1}l_1  & (-1)^{m_2} l_2 & (-1)^{m_3}l_3 & (-1)^{m_4} l_4 \\
 1   &  1      &   1     &   1     \\
\end{array} \right)
\left(\begin{array}{c}
 |u_{1}|^2 \\
       |u_{3}|^2 \\
       |u_{5}|^2\\
       |u_{7}|^2  \\
      \end{array}\right) =
\left(\begin{array}{c}
  0 \\
  0  \\
  0\\
  \frac{1}{2} \\
\end{array}\right).
\end{equation}
This is equivalent to
\begin{equation}         \label{C2 : 60}
\left\{ \begin{array}{l}
2|u_{1}|^2 =\frac{\pm \Delta_{234} }{ \sq}>0,\\
2|u_{3}|^2 =\frac{\pm \Delta_{134}}{ \sq}>0, \\
2|u_{5}|^2 =\frac{\pm \Delta_{124}}{ \sq}>0, \\
2|u_{7}|^2  =\frac{\pm \Delta_{123}}{ \sq}>0,\\
\end{array}\right.
\end{equation}
admitting a unique solution 
in $|u_{1}|^2$, $|u_{3}|^2$, $|u_{5}|^2$, $|u_{7}|^2$. Then, by looking at relations
 (\ref{C1: 26}) we see that
$N(\Theta)$ intersects only one of
 the eigenspaces $\overset{-}{S}{}^{(+,+)}_{(\pm, \pm, \pm)}$.
In particular, in the non empty intersection case,
 dim\, ${S}{}^{(+,+)}_{(\pm, \pm, \pm)}
\cap N(\Theta)= 9$
 and ${S}{}^{(+,+)}_{(\pm, \pm, \pm)}
\cap N(\Theta)/ \tilde{G}$ is diffeomorphic to a $S^2$. $\square$

\vspace{0.5 cm}

We consider now the second mentioned family of singularities on $\mathcal Z^6(\Theta)$, coming from
points $\zw\in \mathbb{H}^8$ having no
 quaternionic pairs in a ${}^{(\pm , \pm)}V_3^{\alpha}$. Observe first that the remaining eigenspaces  ${}^{\pm}V_1^{\alpha}, \,\, {}^{\pm}V_2^{\alpha},\,\, V_4^{\alpha}$ and
 $V_5^{\alpha},\,\, \alpha=1,2,3,4$ are not invariant for the  action
 of the group $\widetilde{G}$. Accordingly, the corresponding strata ${S}\subset$
 $\mathbb{H}^8$ have to be defined as  

 \begin{equation}    \label{C2 : 63}
 \widetilde{G}\cdot  \overbrace{\Bigg\{ \left(\begin{array}{cc|cc|cc|cc}
u_1 & u_2  & u_3 & u_4 & u_5 &  u_6 & u_7 &  u_8\\
\end{array}\right)\Bigg\}}^{V:=},
 \end{equation}

\noindent where $V$ has quaternionic pairs in any of the ${}^{\pm}V_1^{\alpha},$
  ${}^{\pm}V_2^{\alpha},$ $V_4^{\alpha}$ and $V_5^{\alpha} $.   

By examining all cases we get the following:

\begin{prop}
 Let $\zw\in N(\Theta)$ be point in the $V$ defined by
$(\ref{C2 : 63})$.
 Then $\zw$ has at most two
 quaternionic
pairs which are contained  either in  
${}^{\pm}V_1^{\alpha}$ or in ${}^{\pm}V_2^{\alpha}$. 
\end{prop}

As a consequence, we can list the singular strata ${S}$ that will give the 
$22$ spheres mentioned in Theorem B (ii):

\begin{equation}        \label{C2 : 86}
\begin{aligned}
  i)& \quad \quad\,\,\,S{}^{123}_4= \tilde{G}\cdot \Bigg\{
\left( \begin{array}{cc|cc|cc|cc}
z_1 & z_2  & z_3 & z_4   &    z_5 & z_6   & 0    &  0   \\
0   &    0 &  0  & 0     &   0   &  0     & w_7  &  w_8\\
\end{array}\right) \Bigg\}, \\
& \hspace{1 cm} \mathrm{and} \,\,
S{}^{124}_3,
S{}^{134}_2,
S{}^{234}_1,
S{}^1_{234},
S{}^2_{134},
S{}^3_{124},
S{}^4_{123}, \\
 ii)&\, \quad \quad  {}_{\pm}S{}^{123}_4=\tilde{G}\cdot 
\Bigg\{
\left( \begin{array}{cc|cc|cc|cc}
z_1 & z_2  & z_3 & z_4   &    z_5 & z_6    & 0    &  0   \\
0   &    0 &  0  & 0     &   0    &  0     & w_7  & \pm iw_7\\
\end{array}\right) 
\Bigg\}, \\
& \,\, \hspace{1 cm}\mathrm{and} \,\,
{}_{\pm}S{}^{124}_3,
{}_{\pm}S{}^{134}_2,
{}_{\pm}S{}^{234}_1, 
{}^{\pm}S{}^1_{234},
{}^{\pm}S{}^2_{134},  
{}^{\pm}S{}^3_{124},
{}^{\pm}S{}^4_{123}, \\ 
 iii)&\quad \quad \quad  S{}^{12}_{34}= \tilde{G}\cdot \Bigg\{
\left( \begin{array}{cc|cc|cc|cc}
z_1 & z_2  & z_3 & z_4   &    0    &  0    & 0    &  0   \\
0   &    0 &  0  & 0     &   w_5   &  w_6  & w_7  &  w_8\\
\end{array}\right) \Bigg\}, \\
& \hspace{1 cm} \mathrm{and} \,\,
S^{13}_{24},
S^{14}_{23}, 
S^{23}_{14}, 
S^{24}_{13}, 
S^{34}_{12}.\\
 iv)& \quad{}_{(\pm, \pm)}S{}^{12}_{34}=\tilde{G}\cdot 
\Bigg\{
\left( \begin{array}{cc|cc|cc|cc}
z_1 & z_2  & z_3 & z_4   &   0   & 0        & 0    &  0   \\
0   &    0 &  0  & 0     &  w_5  & \pm iw_5 & w_7  & \pm i w_7\\
\end{array}\right) \Bigg\}, \\
&  \,\, \hspace{1 cm} \mathrm{and} \,\,
{}_{(\pm,\pm)}S{}^{13}_{24},
{}_{(\pm,\pm)}S{}^{14}_{23}, 
{}_{(\pm,\pm)}S{}^{23}_{14}, 
{}_{(\pm,\pm)}S{}^{24}_{13},  
{}_{(\pm,\pm)}S{}^{34}_{12}, 
{}^{(\pm,\pm)}S{}^{13}_{24}, \\
& \quad \quad \quad \quad
{}^{(\pm,\pm)}S{}^{14}_{23}, 
{}^{(\pm,\pm)}S{}^{23}_{14}, 
{}^{(\pm,\pm)}S{}^{24}_{13}, 
{}^{(\pm,\pm)}S{}^{34}_{12}. \\
\end{aligned}
\end{equation}

  \begin{lem}
  Let ${S}=\widetilde{G}\cdot V$ be any of
 the strata $(\ref{C2 : 63})$. Then
  $dim\, S= dim\, V + 2$.
  \end{lem}

 $Proof.$ Look at the subgroup $\widetilde{H} = 
Sp(1)\times U(1)^{\epsilon}\times U(1)$
 of $\widetilde{G}$ fixing $V$ and let:

 \begin{equation}    \label{C2 : 63a}
\widetilde{G}\times_{\widetilde{H}} V=\big\{ [(g, \uno)] \,\, | \,\,
 g\in \widetilde{G} \,\, \mathrm{and} \,\,  \uno\in V \big\},
 \end{equation}

 \noindent where  $ [(g, \uno)] =  [(g_1, \uno_1)] $ if  $(g_1, \uno_1)= (gh^{-1},$
   $ h\cdot\uno)$ for some $h\in \widetilde{H}$,
 by definition a $V$-vector bundle over
  $\widetilde{G}/\widetilde{H}$. A $\widetilde{G}-$action
 on $ \widetilde{G}\times \big(\widetilde{G}\times_{\widetilde{H}} V \big)$
  is defined as:
\begin{equation}    \label{C2a : 63b}
 g'\cdot [(g, \uno)]:= [(g'g, \uno)],
 \end{equation}
so that $\widetilde{G}\cdot V\cong
   \widetilde{G}\times_{\widetilde{H}} V$ through the $\widetilde{G}-$equivariant
 diffeomorphism 
\begin{equation}      \label{C2 : 63b}
\begin{aligned}
 \Gamma\, :\widetilde{G}& \cdot V\longrightarrow   \widetilde{G}\times_{\widetilde{H}} V, \\
& g\cdot \uno     \longmapsto   [(g, \uno)].\\
\end{aligned}
\end{equation}
Thus
$\widetilde{G}\cdot V$ can be looked at as a $V$-vector
 bundle over $\widetilde{G}/\widetilde{H} \cong S^2$. $\square$ 
 \newline

\noindent  Next, we have:

 \begin{teor}
 The strata $S{}^{\alpha\beta\gamma}_{\delta}$ and
$S{}^{\alpha\beta}_{\gamma\delta}$
are such that
 \begin{equation}    \label{C2 : 102}
 \begin{aligned}
& i)\,\,\,\,\quad S{}^{\alpha\beta\gamma}_{\delta}\cap N(\Theta) = 
 {}^+S{}^{\alpha\beta\gamma}_{\delta}\cap N(\Theta) \bigcup
 {}^-S{}^{\alpha\beta\gamma}_{\delta}\cap N(\Theta), \\
& ii)\,\,\,\,\quad S{}^{\alpha}_{\beta\gamma\delta}\cap N(\Theta) = 
 {}_+S{}^{\alpha}_{\beta\gamma\delta}\cap N(\Theta) \bigcup
 {}_-S{}^{\alpha}_{\beta\gamma\delta}\cap N(\Theta), \\
& iii)\,\,\,\,\quad  S{}^{\alpha\beta}_{\gamma\delta}\cap N(\Theta)\,\,
   is\,\,connected. \\
 \end{aligned}  \hspace{3 cm}  
 \end{equation}
Moreover
 $\Sigma({\Theta})=\big( \bigcup_{i,j}
 \big(S^{\alpha_i\beta_i\gamma_i}_{\delta_i}\cap
  N(\Theta)\big) \cup
 \big(S^{\alpha_j\beta_j}_{\gamma_j\delta_j}\cap
 N(\Theta)\big)\big)/\widetilde{G}$,
and  each
 $S{}^{\alpha\beta}_{\gamma\delta}$
 contains four substrata
 of those listed in $(\ref{C2 : 86})$, at point $iv)$.
 \end{teor}

$Proof.$
To fix the argument, consider $S^{123}_4$ in (\ref{C2 : 86}) point $i)$.
 The $Sp(1)-$moment map equations and the sphere equation yield a system
\begin{equation}       \label{C2 : 103}
 \left\{\begin{array}{l}
\sum_{\alpha = 1}^6|z_{\alpha}|^2=\frac{1}{2}, \\
\sum_{\alpha=1}^6 (z_{\alpha})^2=0, \\
\end{array} \right. \,\,  \,\, 
 \left\{\begin{array}{l}
|w_7|^2 + |w_8|^2=\frac{1}{2}, \\
(w_7)^2 + (w_8)^2=0, \\
\end{array}\right.
\end{equation}
and by lemma $3.4$, it follows that $dim\, S^{123}_4=18$. This system gives in particular
 $w_8 =\pm
i w_7$ and $S^{123}_4\cap N(\Theta)\subseteq
{}^+S^{123}_4\cup {}^-S^{123}_4$. Note that 
equations
(\ref{C2 : 103}) coincide with the $Sp(1)-$moment map equations
restricted on ${}^{+}S^{123}_4$ and ${}^{-}S^{123}_4$. Also the
$ T^3_{\Theta}-$moment map equations on $S^{123}_4$
\begin{equation}       \label{C2 : 104}
 \left\{\begin{array}{l}
p_1\mathrm{Im}(z_{1}\overline{z}_2) + p_2\mathrm{Im}(z_{3}\overline{z}_4) 
+p_3\mathrm{Im}(z_{5}\overline{z}_6) +p_4\mathrm{Im}(w_{7}\overline{w}_8)=0,\\
 q_1\mathrm{Im}(z_{1}\overline{z}_2)+ q_2\mathrm{Im}(z_{3}\overline{z}_4) 
+q_3\mathrm{Im}(z_{5}\overline{z}_6) +q_4\mathrm{Im}(w_{7}\overline{w}_8)=0, \\
 l_1\mathrm{Im}(z_{1}\overline{z}_2)+ l_2\mathrm{Im}(z_{3}\overline{z}_4) 
+l_3\mathrm{Im}(z_{5}\overline{z}_6) +l_4\mathrm{Im}(w_{7}\overline{w}_8)=0, \\
\end{array}\right. 
\end{equation}
coincide with the ones on ${}^{+}S^{123}_4$ and ${}^{-}S^{123}_4$.
 Since 
$N(\Theta)$ is $\tilde{G}-$ invariant, it follows $S{}^{123}_{4}
\cap N(\Theta)=  {}^+S{}^{123}_{4}\cap
 N(\Theta) \bigcup  
{}^-S{}^{123}_{4}\cap N(\Theta)$. 
If
$ {}^+ \tilde{S}{}^{123}_{4}$
denotes the intersection ${}^+S{}^{123}_{4}\cap
 N(\Theta)$, then:
\begin{equation}
 \left.\begin{array}{c}
 {}^+\tilde{S}{}^{123}_{4} = 
\tilde{G}\cdot
 \Big\{\zw\in \mathbb{H}^8 \big|
\mathrm{Im}(z_{1}\overline{z}_2)= 
\frac{- \Delta_{234} }{4 \Delta{123}},
\mathrm{Im}(z_{3}\overline{z}_4)=
 \frac{ \Delta_{134} }{4 \Delta{123}},
\mathrm{Im}(z_{5}\overline{z}_6)=
 \frac{ -\Delta_{124} }{4 \Delta{123}}, \\
 \sum_{\alpha=1}^6(z_{\alpha})^2=0,
 \sum_{\alpha=1}^6|z_{\alpha}|^2=\frac{1}{2},
 \sum_{\alpha=1}^6|w_{\alpha}|^2=0,
 w_8 = iw_7, 
 |w_7|^2= \frac{1}{4}.
\Big\} 
\end{array}\right.
\end{equation}
\noindent Moreover, $dim\, \big(\,S{}^{123}_4\cap
 N(\Theta)\big)/\tilde{G}=
 dim\,  \big(\,{}^{\pm}S{}^{123}_4\cap 
N(\Theta)\big) /\tilde{G}=2$ and this intersection has two connected components.
Namely, it is easy to see that both of
 $\big( {}^{+} S{}^{123}_4 \cap 
N(\Theta) \big)/ \tilde{G} $ and
 $\big( {}^{-} S{}^{123}_4 \cap 
N(\Theta) \big)/ \tilde{G} $  give twistorial lines
$S^2$. A similar argument applies to any $ S{}^{\alpha\beta\gamma}_{\delta}$ or $ S{}^{\alpha}_{\beta\gamma\delta}$, yielding sixteen twistorial $S^2$. The remaining six $S^2$ come from the strata  $S{}^{\alpha\beta}_{\gamma\delta}$. Refer in particular to
 ${S}{}^{12}_{34}$ defined in (\ref{C2 : 86}), $iii)$.
The $Sp(1)-$moment map 
equations and sphere equation are now
 \begin{equation}         \label{C2 : 107}
   \left\{\begin{array}{l}
\sum_{\alpha=1}^4 |z_{\alpha}|^2 =\frac{1}{2}, \\
\sum_{\alpha=1}^4 (z_{\alpha})^2 =0,       \\
\end{array} \right.\,\,  \,\,
 \left\{\begin{array}{l}       
\sum_{\alpha=5}^8 |w_{\alpha}|^2=\frac{1}{2}, \\
\sum_{\alpha=5}^8(w_{\alpha})^2=0, \\
 \end{array} \right.      \hspace{3 cm} {}
  \end{equation}
and the $T^3_{\Theta}$-moment map equations
 \begin{equation}             \label{C2 : 107a}
 \left\{\begin{array}{l}
p_1\mathrm{Im}(z_{1}\overline{z}_2) + p_2\mathrm{Im}(z_{3}\overline{z}_4) 
+p_3\mathrm{Im}(w_{5}\overline{w}_6) +p_4\mathrm{Im}(w_{7}\overline{w}_8)=0,\\
 q_1\mathrm{Im}(z_{1}\overline{z}_2)+ q_2\mathrm{Im}(z_{3}\overline{z}_4) 
+q_3\mathrm{Im}(w_{5}\overline{w}_6) +q_4\mathrm{Im}(w_{7}\overline{w}_8)=0, \\
 l_1\mathrm{Im}(z_{1}\overline{z}_2)+ l_2\mathrm{Im}(z_{3}\overline{z}_4) 
+l_3\mathrm{Im}(w_{5}\overline{w}_6) +l_4\mathrm{Im}(w_{7}\overline{w}_8)=0. \\
 \end{array}\right. 
\end{equation}

Like in  the previous case we obtain
$\big( S{}^{12}_{34} \cap 
N(\Theta) \big)/ \tilde{G} \cong S^2$. 
Note, that the four strata 
${}_{(\pm,\pm)}S{}^{12}_{34} \subset S{}^{12}_{34} $
and the mentioned argument also gives 
$dim\,{}_{(\pm,\pm)}S{}^{12}_{34} \cap N(\Theta)=7 $ and
 $dim\,\big( {}_{(\pm,\pm)}S{}^{12}_{34}
 \cap N(\Theta)\big)/ \widetilde{G}=0 $ and these give the points listed in statement (iii) of Theorem B in the Introduction. .
 $\square$

\vspace{0.5 cm}

\begin{oss}
There is a real structure on the twistor space
$\mathcal{Z}^6(\Theta)$
coming from the multiplication by the second quaternionic unit $j$
on vectors $\zw\in\mathbb{H}^8$:
\begin{equation}            \label{C2 : 107d}
\mathcal{J}(\zw) = (-\underline{w}, \underline{z}).
\end{equation}
Under this $\mathcal{J}-$map our strata transform according to
$\mathcal{J}(S^{\alpha\beta\gamma}_{\delta})=
\mathcal{J}(S^{\delta}_{\alpha\beta\gamma})$,
$\mathcal{J}({}_{\pm}S^{\alpha\beta\gamma}_{\delta})=
\mathcal{J}({}^{\pm}S^{\delta}_{\alpha\beta\gamma})$,
$\mathcal{J}(S^{\alpha\beta}_{\gamma\delta})=
\mathcal{J}(S^{\gamma\delta}_{\alpha\beta})$ and 
$\mathcal{J}({}_{(\pm,\pm)}S^{\alpha\beta}_{\gamma\delta})=
\mathcal{J}({}^{(\pm,\pm)}S^{\gamma\delta}_{\alpha\beta})$.
\end{oss}

\vspace {1 cm}

\section{The Quotient Orbifolds $\mathcal{O}^4 (\Omega)$
 and their twistor space $\mathcal{Z}^6 (\Omega)$}

Consider now matrices
\begin{align}
& A(\omega_{\alpha})=
  \left( \begin{array}{cc}
  \cos \omega_{\alpha} & \sin \omega_{\alpha}  \\
  -\sin \omega_{\alpha} & \cos \omega_{\alpha}  \\
  \end{array} \right),  \qquad  \qquad
 \omega_{\alpha}=
p_{\alpha} t + q_{\alpha} s,
\end{align} 
with $t, s,  \in [ 0, 2\pi)$, where $\alpha = 1,2,3$ and
\begin{equation}
\Omega =  \left( \begin{array}{cccc}
p_1 & p_2&p_3\\
q_1 & q_2&q_3  
  \end{array} \right) 
\end{equation}
is a matrix of  integral weights.

The  sphere $S^{27}\subset
\mathbb{H}^7$ is acted on by the group
$G^{\Omega }=Sp(1)\times T^2_{\Omega}\subset SO(7)\subset Sp(7)$, 
whose factor $Sp(1)$ acts by left  quaternionic multiplication 
and the second factor $T^2_{\Omega}$ 
by matrices:   

 \begin{equation}      
 A(\Omega)=\left(\begin{array}{c|c|c|c}
1 & 0 & 0 & 0  \\
 \hline
  0 & A(\omega_1) & 0 & 0  \\
  \hline
 0  & 0   & A(\omega_2)  & 0\\
  \hline
 0 & 0 & 0 & A(\omega_3)  \\
 \end{array} \right) \in SO(7),
 \end{equation}

Accordingly, the zero set of the moment map
 $N(\Omega)= 
\mu^{-1}(0)\cap \nu^{-1}(0)\subset S^{27} \subset
\mathbb{H}^7$ 
admits the following $G^{\Omega}-$strata \cite{bg}:
\begin{equation}  \label{C3 : 5}
\begin{aligned}
& S_0 =\{ \uno\in N(\Omega)\, |\, u_1 =0\}, \\
& S_1 =\{ \uno\in N(\Omega)\, |\, u_1 \ne0 \}. \\
\end{aligned}
\end{equation}

The following lemma points out a
minor correction to statements of 
corollary $2.3$ and lemma $3.2$ in \cite{bg}.

\begin{lem}
 The $G^{\Omega}$ action on
$S_1= N(\Omega)\cap \{u_1\ne 0\}$ is
\begin{equation}  \label{C3 : 7}
\begin{aligned}
&  i)\quad  locally \,\, free\,\, if \,\, and\,\,
 only\,\, if \,\,
\Delta_{\alpha\beta}\ne 0, \,\, \forall\, 1\leq \alpha <
 \beta \leq 3,  \\
&  ii) \quad  free\,\, if \,\, and \,\, only\,\, if \,\,
gcd(\Delta_{12}, \Delta_{13}, \Delta_{23})= \pm 1. \hspace{3 cm} {}
\end{aligned}
\end{equation}
(The latter is a weaker condition than the one in
Lemma 3.2 (ii) of \cite{bg}).

\end{lem}

 $Proof.$   
Since  $ u_1\ne 0$, the $Sp(1)$ factor acts trivially, and it is sufficient
to look at the $T^2_{\Omega}-$
action. On the other hand, it is easy to see that, like in Lemma 2.1, no quaternionic pair can vanish on $N(\Omega)$ provided all the minors $\Delta_{\alpha \beta}$ of $\Omega$ are non-zero. 
It follows that the fixed point equations read:
\begin{equation}   \label{C3 : 8}
\overbrace{
\left(\begin{array}{cc}
cos\, \omega_{\alpha} & sin \, \omega_{\alpha} \\
- sin\, \omega_{\alpha} & cos\, \omega_{\alpha} \\
\end{array}  \right)  }^{A(\omega_{\alpha}):=}
\left(\begin{array}{c}
u_{2\alpha} \\
u_{2\alpha+1} \\
\end{array} \right)  =
\left(\begin{array}{c}
u_{2\alpha} \\
u_{2\alpha + 1} \\
\end{array} \right), \quad \alpha = 1,2,3,
\end{equation}
 where $\omega_{\alpha}= p_{\alpha}t + 
q_{\alpha}s$, and $t,s \in [ 0, 2\pi)$, so that $
A(\omega_{\alpha})=id_{2\times 2},$ and
$e^{i(p_{\alpha}t + 
q_{\alpha}s)}=1$,
yielding only the trivial 
solution if and only if
$gcd(\Delta_{12}, \Delta_{13}, 
\Delta_{23})= \pm 1$. 
The locally free conditions were
correctly proved in \cite{bg}. 
$\square$ \newline

Consider now the action of
$\widetilde{G}^{\Omega} = T^2_{\Omega}\times
Sp(1)\times U(1)$ on $N(\Omega)$, similar to the one
in (\ref{C2 : 1}).   
Let $\zw$ be a point in $N(\Omega)
\subset S^{27}\subset \mathbb{H}^7$ 
\begin{equation}  \label{C3 : 23}
  {}^T(\underline{z}, \underline{w} ):=
  \left(\begin{array}{c|cc|cc|cc}
      z_1 & z_2  & z_3 & z_4 & z_5 & z_6 & z_7  \\
      w_1 & w_2 &  w_3 & w_4 & w_5 & w_6 & w_7  \\
    \end{array}\right),
\end{equation}
where $(z_{\beta}, w_{\beta})\in
\mathbb{C}\times\mathbb{C} $, $\beta 
= 1,...,7$. The fixed point
equations for each quaternionic pairs
$(u_{2\alpha}, u_{2\alpha + 1})$ read now: 
\begin{equation}   \label{C3 : 22}
  A(\omega_{\alpha})
  \left(\begin{array}{cc}
      z_{2{\alpha}} & w_{2{\alpha}} \\
      z_{2{\alpha}+1 } & w_{2{\alpha}+1} \\
    \end{array} \right) =
  {\left[ \begin{array}{c}
        \left(\begin{array}{cc}
            \epsilon & -\overline{\sigma}  \\
            \sigma & \overline{\epsilon} \\
          \end{array} \right)
        \left(\begin{array}{cc}
            z_{2\alpha} & z_{2\alpha +1} \\
            w_{2\alpha} & w_{2\alpha +1} \\
          \end{array} \right)
        \left(\begin{array}{cc}
            \rho& 0 \\
            0 & \rho \\
          \end{array} \right)
      \end{array}\right]}^T,
\end{equation}
where $\lambda= \epsilon + 
j\sigma\in Sp(1)$, $\rho\in U(1)$.
These equations can be rewritten 
as:
\begin{equation}  \label{C3 : 24}
  \overbrace{\left( \begin{array}{cc|cc}
        0                          & -\sigma\rho &
        -sin\omega_{\alpha}  & 
        cos\, \omega_{\alpha} - \overline{\epsilon}\rho  \\
-\sigma\rho     & 0                      & 
cos\, \omega_{\alpha} - 
\overline{\epsilon}\rho & sin\omega_{\alpha} \\
\hline
-sin\omega_{\alpha}   & cos\, \omega_{\alpha} - \epsilon\rho
&    0       &      \overline{\sigma}\rho  \\
cos\, \omega_{\alpha} 
-\epsilon \rho & sin\omega_{\alpha}  
&  \overline{\sigma}\rho   &   0    \\
\end{array} \right)}^{\widetilde{M}_{\alpha}:=}
\left(\begin{array}{c}
    z_{2\alpha} \\
    z_{2\alpha +1}\\
    w_{2\alpha} \\
    w_{2\alpha +1} \\
  \end{array} \right) =
\left(\begin{array}{c}
    0 \\
    0\\
    0\\
    0 \\
  \end{array} \right),
\end{equation}
for each $\alpha = 1,2,3$.  Then, similarly to Lemma 3.1 and Proposition 3.1, we get:

\begin{prop}
Let $\zw$ be a point in $N(\Omega)$. Then, 
up to $\widetilde{G}^{\Omega}-$conjugation, we have that
 ${\widetilde{G}^{\Omega}_{\zw}}\subset T^2_{\Omega}
\times \{\lambda\in Sp(1)\,\, |\,\,
 \sigma =0 \}\times U(1)$.
\end{prop}

\begin{prop}  For each $\alpha =1,2,3$,
$det\, \widetilde{M}_{\alpha}=0$ 
if and only if:
\begin{equation} \label{C3 : 63}
\bar \rho e^{\pm i\omega_{\alpha}}= Re( \epsilon ) \pm i\sqrt{ Im( \epsilon)^2 +
|\sigma |^2}.
\end{equation}
\end{prop}

Now. singularities of the twistor space
 $\mathcal{Z}^6(\Omega)$ can be described by looking at the two strata $S_0$ and $S_1$.

\vspace{0.3 cm}

On $S_1 : \{u_1 \ne 0\}$ we have 

\begin{lem}   
  The fixed point equations on $S_1$ and with respect to 
  $u_1$ give:
  \begin{equation}     \label{C3 : 15}
\rho= Re(\epsilon) \pm i\sqrt{Im(\epsilon)^2 + |\sigma|^2}.
\end{equation}
\end{lem}

$Proof.$
In fact:
\begin{equation}    \label{C3 : 16}
u_1 = \lambda u_1 \rho, \,\, \iff \,\,
\left( \begin{array}{cc}
1-\epsilon\rho & \overline{\sigma}\rho   \\
 -\sigma\rho &  1- \overline{\epsilon}\rho  \\
   \end{array} \right)
\left(
 \begin{array}{c}
 z_1 \\
 w_1\\
 \end{array}\right)   =
 \left(
 \begin{array}{c}
 0\\
 0\\
 \end{array}\right).
\end{equation}
 where $\lambda = \epsilon + j\sigma\in Sp(1)$ and $\rho\in U(1)$. There are non trivial solutions if and only if the determinant vanishes, and this gives the stated condition.
$ \square$

\begin{oss}
 If $\sigma = 0$
 in  $(\ref{C3 : 22})$, the fixed point equations with respect to the first quaternionic
  coordinate
 become $\epsilon u_1 \rho = u_1$, i. e.
 $\epsilon z_1 \rho + j\overline{\epsilon}w_1 \rho = z_1 +jw_1.$
Then, when $z_1\ne 0$ and $w_1\ne 0$, we get
 $ \epsilon\rho=1, \,\, and \, \, 
  \epsilon=\rho=\pm 1$, conditions that give
the non effective subgroup.
Moreover, for any
  $(\underline{z}, \underline{w})\in S_1$,
the orbit contains points
with both $z_1 \ne 0$ and  $w_1 \ne 0$.
\end{oss}
 
\noindent Thus, by Proposition 4.1, Lemma 4.2 and the above remark:

 \begin{prop}
The stratum $S_1$ does not give rise to any
 singular point on  $\mathcal{Z}^6 (\Omega)$.
 \end{prop}

\vspace{0.3 cm}

Now, we look at $S_0: \{u_1 =0\}.$ Since the elements of $S_0$ have three quaternionic pairs, this situation can be treated like the zero set $N(\Theta)$ studied in the previous 
section. In fact, by using similar arguments  of
Proposition $3.2$ and Lemma $3.2$, 
we see that the possible solutions for the equations (\ref{C3 : 24}) belong 
to one of the space $ {}^{\pm}V_{1}^{\alpha}$,
$ {}^{\pm}V_{2}^{\alpha}$,
$ {}^{(\pm, \pm)}V_{3}^{\alpha}$,
$ V_{4}^{\alpha}$ or $ V_{5}^{\alpha}$, listed in (\ref{C2 : 30b}).

By using Corollary $3.2$, we can
distinguish two families of strata defined by the action
of $\widetilde{G}^{\Omega}$ on $\mathbb{H}^7$.
The first family is given by
points $\zw$  whose  quaternionc pairs are contained
in ${}^{(\pm ,\pm )}V_3^\alpha$, $\alpha=1,2,3$. Instead,
 points in the second family
have no quaternionic pairs  in ${}^{(\pm ,\pm )}V_3^\alpha$.
The first family of strata has points in one of the following:

 \begin{equation}   \label{C3 : 46}
 \begin{aligned}
& \,\,\, \quad {S}_{(\pm,\pm,\pm)}:= 
\Bigg\{ \left(\begin{array}{c|cc|cc|cc}
0  &  z_3 &  \pm iz_3 & z_5 &  \pm iz_5 & z_7 &  \pm iz_7\\
0  &  w_3 &  \pm iw_3 & w_5 &  \pm iw_5 & w_7 &  \pm iw_7 \\
\end{array}\right)\Bigg\}, \\
\end{aligned}
\end{equation}
The second family gives rise to:
\begin{equation}   \label{C3 : 47}
 \begin{aligned}
& i) \,\,\,  \quad S{}^{12}_3:=\widetilde{G}^{\Omega}\cdot \Bigg\{
\left( \begin{array}{c|cc|cc|cc}
0   & z_2  & z_3 & z_4  &    z_5 &  0    &  0   \\
0   &    0 &  0  & 0    &   0    & w_6  &  w_7\\
\end{array}\right) \Bigg\},   \\
& ii) \,\,\, \quad  {}_{\pm}S{}^{12}_3:=\widetilde{G}^{\Omega}\cdot \Bigg\{
\left( \begin{array}{c|cc|cc|cc}
0   & z_2  & z_3 & z_4   &    z_5 &  0    &  0   \\
0   &    0 &  0  & 0     &   0    & w_6   &  \pm iw_6\\
\end{array}\right) \Bigg\},     \\
& \mathrm{and} \,
S{}^{13}_2, S{}^{23}_1, S{}^1_{23},S{}^2_{13},S{}^3_{12}, \\
& \mathrm{as\; well \; as} \,
{}_{\pm}S{}^{13}_2, {}_{\pm}S{}^{23}_1, {}_{\pm}S{}^1_{23},
{}_{\pm}S{}^2_{13},{}_{\pm}S{}^3_{12},{}{\pm}S{}^{12}_3, 
 {}^{\pm}S{}^{13}_2, {}^{\pm}S{}^{23}_1, {}^{\pm}S{}^1_{23},
{}^{\pm}S{}^2_{13},{}^{\pm}S{}^3_{12}. \\
\end{aligned}
\end{equation}
  
We now see the 2-sphere appearing in Theorem C (i).
In fact, through the same argument used in Proposition 3.6 we get:

\begin{prop}
  Just  one of the strata $S_{(\pm,\pm,\pm)}$ intersects
  the zero set
  $N(\Omega)$. This intersection generates a $2-$sphere $S^2$ 
on the twistor space $\mathcal{Z}^6 (\Omega)$
\end{prop}

\noindent As for
the second family of singular strata: 

 \begin{teor}
   The strata
   listed in $(\ref{C3 : 47})$ have no empty
   intersection with the submanifold $N(\Omega)$. Moreover:
   \begin{equation}    \label{C3 : 48}
     \begin{aligned}
       & i)\,\,\,\,\quad S{}^{\alpha\beta}_{\gamma}\cap N(\Omega) = 
       {}^+S{}^{\alpha\beta}_{\gamma}\cap N(\Omega) \bigcup
       {}^-S{}^{\alpha\beta}_{\gamma}\cap N(\Omega), \\
       & ii)\,\,\,\,\quad S{}^{\alpha}_{\beta\gamma}\cap N(\Omega) = 
       {}_+S{}^{\alpha}{\beta\gamma}\cap N(\Omega) \bigcup
       {}_-S{}^{\alpha}{\beta\gamma}\cap N(\Omega), \\
\end{aligned}  \hspace{3 cm}  
\end{equation}
and each of the connected components ${}^{\pm}S{}^{\alpha\beta}_{\gamma}\cap N(\Omega)$
and ${}_{\pm}S{}^{\alpha}_{\beta\gamma}\cap N(\Omega)$
generate a singular point on the twistor space  $\mathcal{Z}^6 (\Omega)$.
\end{teor}

The proof is a consequence of the same dimensional argument
 used in Theorem $3.1$. The connected components in (\ref{C3 : 48}) $i)$ and $ii)$,
 provide the singular points' description in 
Theorem C (ii).

\end{document}